%% file: main.tex
\documentclass[a4paper, 11pt, reqno]{amsart}
\usepackage[utf8]{inputenc}
\usepackage[english]{babel}
\usepackage{mathtools}
\usepackage{hyperref}
\usepackage{subfiles}
\usepackage{amsthm, amsmath, amssymb, amsfonts, graphicx}
\hypersetup{colorlinks=true,urlcolor=blue,citecolor=red}

\newtheorem{theorem}{Theorem}[section]

\newtheorem{corollary}[theorem]{Corollary}

\newtheorem{definition}[theorem]{Definition}
\newtheorem{example}[theorem]{Example}

\newtheorem{lemma}[theorem]{Lemma}

\newtheorem{remark}[theorem]{Remark}

\numberwithin{equation}{section}

\title{On Collatz Conjecture}
\author{Erhan TEZCAN}
\date{March 12, 2019}
\email{erhany96@gmail.com, etezcan19@ku.edu.tr}

\begin{document}

\begin{abstract}
The Collatz Conjecture can be stated as: using the reduced Collatz function $C(n) = (3n+1)/2^x$ where $2^x$ is the largest power of 2 that divides $3n+1$, any odd integer $n$ will eventually reach 1 in $j$ iterations such that $C^j(n) = 1$. In this paper we use reduced Collatz function and reverse reduced Collatz function. We present odd numbers as sum of fractions, which we call ``fractional sum notation'' and its generalized form ``intermediate fractional sum notation'', which we use to present a formula to obtain numbers with greater Collatz sequence lengths. We give a formula to obtain numbers with sequence length 2. We show that if trajectory of $n$ is looping and there is an odd number $m$ such that $C^j(m) = 1$, $n$ must be in form $3^j\times2k + 1, k \in \mathbb{N}_0$ where $C^j(n) = n$. We use Intermediate fractional sum notation to show a simpler proof that there are no loops with length 2 other than trivial cycle looping twice. We then work with reverse reduced Collatz function, and present a modified version of it which enables us to determine the result in modulo 6. We present a procedure to generate a Collatz graph using reverse reduced Collatz functions.
\bigskip\par \textbf{Keywords: } Number Theory, 3x+1 Problem, Collatz Conjecture, Unsolved Problem, Sequence, Cycle, Loop, Reverse, Reduced Collatz
\par \textbf{MSC:} 11B37 (Primary); 39B05 (Secondary)
\end{abstract}

\maketitle

\subfile{Sections/1-Introduction}

\subfile{Sections/1.1-Reformulation}

\subfile{Sections/2-FSN}

\subfile{Sections/2.1-OnePiece}

\subfile{Sections/2.2-TwoPiece}

\subfile{Sections/2.3-AdditiveGeneration}

\subfile{Sections/3-IFSN}

\subfile{Sections/3.1-Loops}

\subfile{Sections/3.2-Divergence}

\subfile{Sections/4-ReverseFormula}

\subfile{Sections/4.1-Remainders}

\subfile{Sections/4.2-OneTree}

\subfile{Sections/5-Conclusion}

\bibliographystyle{plain}
\bibliography{main}

\end{document}

%% file: Sections/1-Introduction.tex
\section{Introduction}
The Collatz Conjecture is a conjecture, also known as Ulam Conjecture, Syracuse Problem, Hasse's Algorithm and possibly the most descriptively $3x+1$ Problem among many other names. The Collatz function $F : \mathbb{N} \xrightarrow{} \mathbb{N}$ is defined by
\begin{equation}
\label{EQ: Definition of Collatz Function F(n)}
    F(n) = 
    \left\{ 
    \begin{array}{rl}
        n/2     &\mbox{ if $x \equiv 0 \pmod{2}$} \\
        3n+1    &\mbox{ if $x \equiv 1 \pmod{2}$}
    \end{array} 
    \right.
\end{equation}

Let $F^i(n)$ denote the result of $i$ iterations of function $F$ on $x$. It is conjectured that there exists a number of iterations $j$ such that $F^j(n) = 1$. Whether this is true or not is an unsolved problem \cite{ianstewart}, it still is under spotlight when it comes to simple to express yet hard to solve problems. A bibliography regarding Collatz Conjecture can be read in this excellent paper by J. C. Lagarias \cite{lagarias-bib}. It has been computed that every number up to $2^{60}$ do end up reaching $1$ \cite{ericroosendall}. Notice how after $F^j(n) = 1$ the next iteration is $F^{j+1}(n) = 3(1)+1 = 4$, $F^{j+2}(n) = 4/2 = 2$ and $F^{j+3}(n) = 2/2 = 1$. Upon reaching $1$ at some iteration, the remaining iterations will loop on these numbers: $1\xrightarrow{}4\xrightarrow{}2\xrightarrow{}1\xrightarrow{}4\xrightarrow{}2\xrightarrow{}1\xrightarrow{}\ldots$ This loop is often referred to as ``trivial cycle'' \cite{lagarias-bib}. When we refer to $j^{th}$ iteration such that $F^j(n) = 1$ we prefer to assume that this iteration is the first time the result yields 1. 
\begin{corollary}
Let $n, i, j \in  \mathbb{N}$ and $F^j(n) = 1$ where $F^i(n) \ne 1, 1\leq i<j$. We can say that for any $m \in  \mathbb{N}$ such that $m = F^i(n), i<j$ the iteration $F^k(m) = 1$ where $k = j - i$.
\end{corollary}
This corollary formally states that if a number $n$ reaches $1$ in $j$ iterations, every number in that trajectory also reaches $1$.
\begin{corollary}
\label{COROL: Reverse Iteration}
Let $n, i, j \in  \mathbb{N}$ and $F^j(n) = 1$ where $F^i(n) \ne 1, 1\leq i<j$. Then $F^{j+1}(2n) = 1$ for every $n$ and $F^{j+1}((n-1)/3)$ where $n - 1 \equiv 0 \pmod{3}$. This reverse operation will be a point of discussion in the later stages of our paper.
\end{corollary}
\begin{definition}
\label{DEF: Trajectory and Sequence}
There are two phrases we will be using most often: Collatz trajectory and Collatz sequence. A Collatz trajectory is a sequence of numbers starting at $n \in  \mathbb{N}$ and ends at $m \in  \mathbb{N}$ where $m = F^i(n), i \in  \mathbb{N}$. A Collatz Sequence is a sequence of numbers starting from any number $n \in \mathbb{N}$ and ending at $1$, so it is the special case of the Collatz trajectory.   
\end{definition}
We can show the Collatz trajectory using numbers and arrows where each arrow points to an iteration over $F$. 
$$
n \xrightarrow{} F^1(n) \xrightarrow{} F^2(n) \xrightarrow{} \ldots \xrightarrow{} F^{i-1}(n) \xrightarrow{} m
$$
$n$ and $m$ are as defined in definition \ref{DEF: Trajectory and Sequence}. Similarly, a Collatz sequence can be shown as
$$
n \xrightarrow{} F^1(n) \xrightarrow{} F^2(n) \xrightarrow{} \ldots \xrightarrow{} F^{j-1}(n) \xrightarrow{} 1
$$
\begin{example}
The Collatz sequence for $n=3$ is shown as
$$
3 \xrightarrow{} 10 \xrightarrow{} 5\xrightarrow{} 16\xrightarrow{} 8\xrightarrow{} 4\xrightarrow{} 2\xrightarrow{} 1
$$
\end{example}

%% file: Sections/1.1-Reformulation.tex
\subsection{Reformulation of the problem}
A widely used reformulation is to omit the even numbers from the trajectories using a different function. If the conjecture is proved for odd numbers it is intrinsically proved for even numbers, as also seen in corollary \ref{COROL: Reverse Iteration}. To do this we use the reduced Collatz function $C : \mathbb{O}_+ \xrightarrow{} \mathbb{O}_+$ defined by
\begin{equation}
    C(n) = \frac{3n+1}{2^a}
\end{equation}
where $2^a \mid 3n+1$ for the largest value of $a \in \mathbb{N}$ and $\mathbb{O}_+$ is the set of positive odd integers. The notation $m \mid x$ states ``$m$ divides $x$'' or as a congruence $x \equiv 0 \pmod{m}$, for the opposite purpose $m \nmid x$ states that ``$m$ does not divide $x$''. Letter $C$ is preferred with respect to R. E. Crandall \cite{crandall} who most notably used this reformulation, plus it is the first letter of Collatz! We will define a reduced Collatz trajectory where we will be able to see the exponent $a$ explicitly. 
\begin{definition}
The reduced Collatz trajectory using numbers and arrows where each arrow points to an iteration over $C$. 
$$
n \xrightarrow{a_1} C^1(n) \xrightarrow{a_2} C^2(n) \xrightarrow{a_3} \ldots \xrightarrow{a_{i-1}} C^{i-1}(n) \xrightarrow{a_i} m
$$
The reduced Collatz sequence is again similar,
$$
n \xrightarrow{a_1} C^1(n) \xrightarrow{a_2} C^2(n) \xrightarrow{a_3} \ldots \xrightarrow{a_{j-1}} C^{j-1}(n) \xrightarrow{a_j} 1
$$
We say that $n$ has a ``reduced Collatz trajectory length'' $j$ for $C^j(n) = m$ in the trajectory written above. The reduced Collatz sequence has it $m=1$, so for that we prefer to say $n$ has a reduced Collatz sequence length $j$. In short, if we use the word ``sequence'' it describes that the end of the trajectory is $1$.
\end{definition}

\begin{example}
The reduced Collatz sequence for $n=3$ is shown as
$$
3 \xrightarrow{1}  5\xrightarrow{4} 1
$$
The reduced Collatz sequence length of $3$ is $2$.
\end{example}

%% file: Sections/2-FSN.tex
\section{Fractional Sum Notation}
We will present a way to represent numbers that eventually result in $1$ iterating over $C$.
\begin{definition}
\label{DEF: FSN Definition}
Let $n \in \mathbb{O}_+$, $C^j(n) = 1$. The reduced Collatz trajectory is
$$
n \xrightarrow{a_1} C^1(n) \xrightarrow{a_2} C^2(n) \xrightarrow{a_3}\ldots\xrightarrow{a_{j-1}}C^{j-1}(n) \xrightarrow{a_j}1
$$
Then we can write $n$ as
\begin{equation}
\label{DEF: FSN}
    n = \frac{2^{a_1+a_2+\ldots+a_j}}{3^j}-\frac{2^{a_1+a_2+\ldots+a_{j-1}}}{3^j}-\frac{2^{a_1+a_2+\ldots+a_{j-2}}}{3^{j-1}}-\ldots-\frac{2^{a_1}}{3^2}-\frac{2^{0}}{3^1}
\end{equation}
We call this notation as the fractional sum notation, though it has been noted as an ``inverse canonical form'' by G. Helms \cite{ghelms}. Note that a number $n$ with reduced Collatz sequence length $j$ will have $j+1$ terms as seen in the equation above.
\end{definition}
We should explain the terms on the sum of fractions in the equation \eqref{DEF: FSN}. First of all, for every term, the exponent of its numerator is greater than the exponent of the numerator of the term to its right. This is clearly seen by looking at the sums at the exponents, the rightmost sum being 0 as well. The rightmost term is very important: if the exponent of its numerator is greater than 0, the whole expression is an even number. In fact, we defined $n$ to be an odd number, but an even number $m = 2^x \times n$ can be represented similarly
\begin{equation*}
m = \frac{2^{a_1+a_2+\ldots+a_j+x}}{3^j}-\frac{2^{a_1+a_2+\ldots+a_{j-1}+x}}{3^j}-\frac{2^{a_1+a_2+\ldots+a_{j-2}+x}}{3^{j-1}}-\ldots-\frac{2^{a_1+x}}{3^2}-\frac{2^{x}}{3^1}
\end{equation*}
If the rightmost term have exponent 0 then $n$ is an odd number. The first iteration $C(n)$ with exponent $a_1$ will do the operation $(3n+1)/2^{a_1}$. What we get is 
\begin{equation*}
    C(n) = \frac{3n+1}{2^{a_1}}= \frac{2^{a_2+\ldots+a_j}}{3^{j-1}}-\frac{2^{a_2+\ldots+a_{j-1}}}{3^{j-1}}-\frac{2^{a_2+\ldots+a_{j-2}}}{3^{j-2}}-\ldots-\frac{2^{0}}{3^1}
\end{equation*}
In other words, this representation is an explicit expression that shows all the information about iterations of $n$ over $C$. The exponent of denominator of the left-most term tells us how many $3n+1$ operations would happen if $n$ were to iterate over $F$. The exponents of numerator also show how many divisions by 2 occur at a certain iteration. This notation therefore gives us the ability to use the Collatz trajectory all together in a mathematical expression.
\begin{example}
The reduced Collatz sequence for $n=3$ is shown as $3 \xrightarrow{1}  5\xrightarrow{4} 1$. The fractional sum notation for $3$ is then
$$
3 = \frac{2^5}{3^2}-\frac{2^1}{3^2}-\frac{2^0}{3^1}
$$
\end{example}

%% file: Sections/2.1-OnePiece.tex
\subsection{Odd Numbers with reduced Collatz Sequence length 1}
Consider a number $n$ where $C(n) = 1$, the reduced Collatz sequence being $n \xrightarrow{a_1} 1$. The fractional sum notation is then
\begin{equation}
    n = \frac{2^{a_1}}{3^1} - \frac{2^0}{3^1}
\end{equation}
This reduces to
\begin{equation}
\label{EQ: FSN for length 1}
    n = \frac{2^{a_1}-1}{3}
\end{equation}
Since $n$ is an odd integer, equation \eqref{EQ: FSN for length 1} tells us that
\begin{equation*}
    2^{a_1} - 1 \equiv 0 \pmod{3}
\end{equation*}
\begin{remark}
Notice that in equation \eqref{EQ: FSN for length 1} if we try to get $1$ at the right-hand side of the equation the left-hand side becomes $C(n)$.
\begin{equation*}
    \frac{3n+1}{2^{a_1}} = C(n) = 1
\end{equation*}
\end{remark}
\begin{lemma}
\label{LEMMA: Length 1 Odd power}
    $2^{2k+1} - 1\equiv 1 \pmod{3}$, $k \in \mathbb{N}_0$\footnote{$\mathbb{N}_0 = \mathbb{N} \cup \{0\}$}.
\end{lemma}
\begin{proof}
For $k=0$ we have $2^{3} - 1 \equiv 1 \pmod{3}$ which is correct. We assume that the congruence holds for $k=n$ therefore $2^{2n+1} - 1\equiv 1 \pmod{3}$. Looking at $k=n+1$ we have $2^{2(n+1)+1} - 1\equiv 1 \pmod{3}$.
\begin{flalign*}
    2^{2n+2+1}-1&=2^{2n+1}\times4-1 \\
    &=2^{2n+1}\times(3+1)-1 \\
    &=3\times2^{2n+1}+2^{2n+1}-1
\end{flalign*}
Clearly $3\times2^{2n+1}+2^{2n+1}-1 \equiv 1 \pmod{3}$ is correct by induction.
\end{proof}
\begin{lemma}
\label{LEMMA: Length 1 Even power}
    $2^{2k+2} - 1\equiv 0 \pmod{3}$, $k \in \mathbb{N}_0$.
\end{lemma}
\begin{proof}
For $k=0$ we have $2^{2} - 1\equiv 0 \pmod{3}$ which is correct. We assume that the congruence holds for $k=n$ therefore $2^{2n+2} - 1\equiv 0 \pmod{3}$. Looking at $k=n+1$ we have $2^{2(n+1)+2} - 1\equiv 0 \pmod{3}$.
\begin{flalign*}
    2^{2n+2+2}-1&=2^{2n+2}\times4-1 \\
    &=2^{2n+2}\times(3+1)-1 \\
    &=3\times2^{2n+2}+2^{2n+2}-1
\end{flalign*}
Clearly $3\times2^{2n+2}+2^{2n+2}-1 \equiv 1 \pmod{3}$ is correct by induction.
\end{proof}
\begin{theorem}
\label{THEO: Formula for Length 1 theorem}
$C(n) = 1$ where 
\begin{equation}
\label{EQ: Formula for Length 1}
    n = \frac{2^{2k+2}}{3^1} - \frac{2^0}{3^1}, k \in \mathbb{N}_0
\end{equation}
In other words, positive odd integers in the form above have reduced Collatz sequence length 1.
\end{theorem}
\begin{proof}
According to equation \eqref{EQ: FSN for length 1} and lemmas \ref{LEMMA: Length 1 Odd power} and \ref{LEMMA: Length 1 Even power} this is correct. 
\end{proof}
It is important to mention that R. P. Steiner \cite{steiner} proved that the trivial cycle $1\xrightarrow{2}1$ is the only loop with length 1.

%% file: Sections/2.2-TwoPiece.tex
\subsection{Odd Numbers with reduced Collatz Sequence length 2}
Consider a number $n$ where $C^2(n) = 1$, the reduced Collatz sequence being $n \xrightarrow{a_1} C(n) \xrightarrow{a_2} 1$. The fractional sum notation is then
\begin{equation}
\label{EQ: FSN for length 2}
    n = \frac{2^{a_1+a_2}}{3^2}-\frac{2^{a_1}}{3^2}-\frac{2^0}{3^1}
\end{equation}
The values of $a_1$ and $a_2$ for the first 22 odd numbers with reduced Collatz sequence length 2 are given in table \ref{TABLE: Odd Numbers with reduced Collatz Sequence length 2}. We see that for even values of $a_1$ the first $a_2$ for each even value of $a_1$ is $8$, then $14$ which is $8+6$. Similarly, looking at odd values of $a_1$ the first $a_2$ for each odd value of $a_1$ is $4$, then $10$ which is $4+6$.
\begin{table}[!htbp]
\caption{Table of the first 22 odd numbers with Collatz sequence length 2}\label{TABLE: Odd Numbers with reduced Collatz Sequence length 2}
\renewcommand\arraystretch{1.5}
\noindent\[
\begin{array}{|cccc|}
\hline
{n}&{a_1}&{a_2}&{a_1 \text{ Parity}}\\
\hline \hline
 3 & 1 & 4 & \text{Odd} \\ 
 \hline
 13 & 3 & 4 & \text{Odd} \\
 \hline
 53 & 5 & 4 & \text{Odd} \\
 \hline
 113 & 2 & 8 & \text{Even} \\
 \hline
 213 & 7 & 4 & \text{Odd} \\ 
 \hline
 227 & 1 & 10 & \text{Odd} \\ 
 \hline
 453 & 4 & 8 & \text{Even} \\ 
 \hline
 853 & 9 & 4 & \text{Odd} \\ 
 \hline
 909 & 3 & 10 & \text{Odd} \\ 
 \hline
 1813 & 6 & 8 & \text{Even} \\ 
 \hline
 3413 & 11 & 4 & \text{Odd} \\
 \hline
 3637 & 5 & 10 & \text{Odd} \\
 \hline
 7253 & 8 & 16 & \text{Even} \\
 \hline
 7281 & 2 & 14 & \text{Even} \\
 \hline
 13653 & 13 & 4 & \text{Odd} \\
 \hline
 14549 & 7 & 10 & \text{Odd} \\
 \hline
 14563 & 1 & 16 & \text{Odd} \\
 \hline
 29013 & 10 & 8 & \text{Even} \\
 \hline
 29125 & 4 & 14 & \text{Even} \\
 \hline
 54613 & 15 & 4 & \text{Odd} \\
 \hline
 58197 & 9 & 10 & \text{Odd} \\
 \hline
 58253 & 3 & 16 & \text{Odd} \\
 \hline
\end{array}
\]
\end{table}
\begin{theorem}
\label{THEO: Formula for Length 2 theorem}
$C^2(n) = 1$ where
\begin{equation}
\label{EQ: Formula for Length 2}
    n = \frac{2^{a_1+6b+2(-1)^{a_1 \bmod{2}}}}{3^2}-\frac{2^{a_1}}{3^2}-\frac{2^0}{3^1}
\end{equation}
and $a_1, b \in \mathbb{N}$.
\end{theorem}
\begin{proof}
Equation \eqref{EQ: Formula for Length 2} tells us that the congruence 
\begin{equation}
\label{EQ:  Formula for Length 2 Congruence}
2^{a_1+6b+2(-1)^{a_1 \bmod{2}}}-2^{a_1}-3 \equiv 0 \pmod{3^2}, a_1 \in \mathbb{N}, b \in \mathbb{N}_0
\end{equation}
must hold. To show this we will be considering two separate equations, one for the even values of $a_1$ and the other for odd values of $a_1$. This will enable us to use the lemmas \ref{LEMMA: Length 1 Odd power} and \ref{LEMMA: Length 1 Even power} which we proved already. If $a_1$ is even then it is in the form $2k+2, k \in \mathbb{N}_0$. The congruence \eqref{EQ:  Formula for Length 2 Congruence} becomes
\begin{equation*}
    2^{2k+2+6b+2}-2^{2k+2}-3 \equiv 0 \pmod{3^2}
\end{equation*}
Looking at $b=1$ we have $2^{2k+2+6+2}-2^{2k+2}-3 \equiv 0 \pmod{3^2}$. This expression reduces to
\begin{flalign*}
    2^{2k+2+6+2}-2^{2k+2}-3&= 2^{2k+2}\times(2^8-1)-3 \\
    &=2^{2k+2}\times255-3 \\
    &= 2^{2k+2}\times85\times3-3 \\
    &= 3(2^{2k+2}\times85-1) \\
    &= 3(2^{2k+2}\times(84+1)-1) \\
    &= 3\times84\times2^{2k+2}+3(2^{2k+2}-1) \\
    &= 9\times28\times2^{2k+2}+3(2^{2k+2}-1)
\end{flalign*}
The first term is divisible by $3^2$. As for the second term, lemma \ref{LEMMA: Length 1 Even power} shows us $2^{2k+2}-1 \equiv 0 \pmod{3}$ therefore $3(2^{2k+2}-1) \equiv 0 \pmod{3^2}$. We assume that the congruence holds for $b=n$ therefore $2^{2k+2+6n+2}-2^{2k+2}-3 \equiv 0 \pmod{3^2}$. Looking at $b=n+1$ we have $2^{2k+2+6(n+1)+2}-2^{2k+2}-3 \equiv 0 \pmod{3^2}$. This expression reduces to
\begin{flalign*}
    2^{2k+2+6(n+1)+2}-2^{2k+2}-3&=2^{2k+2}(2^{6n+6+2}-1)-3 \\
    &=2^{2k+2}(2^{6n+2}\times2^6-1)-3 \\
    &=2^{2k+2}(2^{6n+2}\times(2^6-1+1)-1)-3 \\
    &=2^{2k+2}(2^{6n+2}\times63+2^{6n+2}-1)-3 \\
    &=2^{2k+2}(2^{6n+2}\times9\times7)+2^{2k+2}\times(2^{6n+2}-1)-3 \\
    &=2^{2k+2}(2^{6n+2}\times9\times7)+2^{2k+2+6n+2}-2^{2k+2}-3
\end{flalign*}
The first term is divisible by $3^2$, rest of the terms are exactly the same as $b=n$ case which by induction proves that the congruence \eqref{EQ:  Formula for Length 2 Congruence} holds for even values of $a_1$. If $a_1$ is odd then it is in the form $2k+1, k \in \mathbb{N}_0$. The congruence \eqref{EQ:  Formula for Length 2 Congruence} becomes
\begin{equation*}
    2^{2k+1+6b-2}-2^{2k+1}-3\equiv 0 \pmod{3^2}
\end{equation*}
Looking at $b=1$ we have $2^{2k+1+6-2}-2^{2k+1}-3 \equiv 0 \pmod{3^2}$. This expression reduces to
\begin{flalign*}
    2^{2k+1+6-2}-2^{2k+1}-3&=2^{2k+1}\times(2^{4}-1)-3 \\
    &=2^{2k+1}\times15-3 \\
    &=2^{2k+1}\times(3\times5)-3 \\
    &=3\times(2^{2k+1}\times(5)-1) \\
    &=3\times(2^{2k}\times(10)-1) \\
    &=3\times(2^{2k}\times(9+1)-1) \\
    &=9\times3\times2^{2k}+3\times(2^{2k}-1)
\end{flalign*}
The first term is divisible by $3^2$. As for the second term, if $k=0$ then $3(2^{0}-1) = 0$, if $k>1$ we know that $2^{2k}-1 \equiv 0 \pmod{3}$ from lemma \ref{LEMMA: Length 1 Even power} and therefore $3(2^{2k}-1) \equiv 0 \pmod{3^2}$. We assume that the congruence holds for $b=n$ therefore $2^{2k+1+6n-2}-2^{2k+1}-3 \equiv 0 \pmod{3^2}$. Looking at $b=n+1$ we have $2^{2k+1+6(n+1)-2}-2^{2k+1}-3 \equiv 0 \pmod{3^2}$. This expression reduces to
\begin{flalign*}
    2^{2k+1+6(n+1)-2}-2^{2k+1}-3&=2^{2k+1}\times(2^{6n-2+6}-1)-3 \\
    &=2^{2k+1}\times(2^{6n-2}\times2^6-1)-3 \\
    &=2^{2k+1}\times(2^{6n-2}\times(2^6+1-1)-1)-3 \\
    &=2^{2k+1}\times(2^{6n-2}\times(2^6-1)+2^{6n-2}-1)-3 \\
    &=2^{2k+1}\times2^{6n-2}\times63+2^{2k+1}\times(2^{6n-2}-1)-3 \\
    &=2^{2k+1}\times2^{6n-2}\times9\times7+2^{2k+1+6n-2}-2^{2k+1}-3
\end{flalign*}
The first term is multiplied by 9 so it is divisible by $3^2$, rest of the terms are exactly the same as $b=n$ case which by induction proves that the congruence \eqref{EQ:  Formula for Length 2 Congruence} holds for odd values of $a_1$. The proofs for even and odd values of $a_1$ therefore prove this theorem.
\end{proof}
\begin{remark}
A pattern similar to what we briefly showed in paragraph below equation \ref{EQ: FSN for length 2} is subtly noticed for odd numbers with Collatz sequence length 3 too. Just like we did in table \ref{TABLE: Odd Numbers with reduced Collatz Sequence length 2}, sorting the numbers in ascending order and look at the exponents $a_1, a_2$ and $a_3$, one can see a similar yet more convoluted pattern. We were unable to generalize this.
\end{remark}

%% file: Sections/2.3-AdditiveGeneration.tex
\section{Using fractional sum notation to obtain numbers with higher sequence lengths}
In this section we will see how we can use the fractional sum notation to obtain an odd number with higher sequence length. Consider a number $n$ with the reduced Collatz sequence $n\xrightarrow{a_1}C(n)\xrightarrow{a_2}C^2(n)\ldots\xrightarrow{a_{j-1}}C^{j-1}(n)\xrightarrow{a_j}1$ of length $j$. We will present a formula to obtain an odd number $m$ with the reduced Collatz sequence $m\xrightarrow{a_1}C(m)\xrightarrow{a_2}C^2(m)\ldots\xrightarrow{a_{j-1}}C^{j-1}(m)\xrightarrow{a_j}C^j(m)\xrightarrow{a_{j+1}}1$ of length $j+1$. This is different than obtaining a number by going reverse, as shown in corollary \ref{COROL: Reverse Iteration}. Consider the same $n$ before but now suppose $n = 3k+1$, in other words $C(k) = n$. Let the exponent be $a_k$ for that iteration, then the reduced Collatz sequence for $k$ is $k\xrightarrow{a_k}n\xrightarrow{a_1}C(n)\xrightarrow{a_2}C^2(n)\ldots\xrightarrow{a_{j-1}}C^{j-1}(n)\xrightarrow{a_j}1$. Notice how the new exponent is at the end of the trajectory for $m$ and at the beginning of the trajectory for $k$.

\begin{theorem}
\label{THEO: Additive Formula}
Let $n, m \in \mathbb{N}$ where $n$ has reduced Collatz sequence length $j$ and $m$ has reduced Collatz sequence length $j+1$.
\begin{equation}
\label{EQ: Additive Formula}
    m = \frac{2^{a_1+a_2+\ldots+a_j}}{3^{j+1}}\times(2^{a_{j+1}}-2^2)+n
\end{equation}
where $a_{j+1} = 2b3^j + 2, b \in \mathbb{N}_0$.
\end{theorem}
\begin{proof}
As defined in definition \ref{DEF: FSN Definition} we have 
\begin{equation}
    n = \frac{2^{a_1+a_2+\ldots+a_j}}{3^j}-\frac{2^{a_1+a_2+\ldots+a_{j-1}}}{3^j}-\frac{2^{a_1+a_2+\ldots+a_{j-2}}}{3^{j-1}}-\ldots-\frac{2^{a_1}}{3^2}-\frac{2^{0}}{3^1}
\end{equation}
By the same definition $m$ is
\begin{equation}
    m = \frac{2^{a_1+a_2+\ldots+a_{j+1}}}{3^{j+1}}-\frac{2^{a_1+a_2+\ldots+a_j}}{3^{j+1}}-\frac{2^{a_1+a_2+\ldots+a_{j-1}}}{3^j}-\ldots-\frac{2^{a_1}}{3^2}-\frac{2^{0}}{3^1}
\end{equation}
Notice that the exponents $a_1, a_2, \ldots, a_j$ are mutual. Now suppose $x$ and $y$ are
\begin{equation*}
    x = \frac{2^{a_1+a_2+\ldots+a_j}}{3^j}
\end{equation*}
\begin{equation*}
    y = -\frac{2^{a_1+a_2+\ldots+a_{j-1}}}{3^j}-\frac{2^{a_1+a_2+\ldots+a_{j-2}}}{3^{j-1}}-\ldots-\frac{2^{a_1}}{3^2}-\frac{2^{0}}{3^1}
\end{equation*}
Clearly $n = x + y$.
\begin{equation*}
    m = \frac{2^{a_1+a_2+\ldots+a_{j+1}}}{3^{j+1}}-\frac{x}{3}+y
\end{equation*}
Adding $x-x$ to the right side gives
\begin{equation*}
    m = \frac{2^{a_1+a_2+\ldots+a_{j+1}}}{3^{j+1}}-\frac{x}{3}-x+x+y
\end{equation*}
Since $n = x + y$
\begin{equation*}
    m = \frac{2^{a_1+a_2+\ldots+a_{j+1}}}{3^{j+1}}-\frac{4x}{3}+n
\end{equation*}
Writing $x$ explicitly we have
\begin{equation*}
    m =\frac{2^{a_1+a_2+\ldots+a_{j+1}}}{3^{j+1}}-\frac{4\times \frac{2^{a_1+a_2+\ldots+a_j}}{3^j}}{3}+n
\end{equation*}
\begin{equation*}
    m =\frac{2^{a_1+a_2+\ldots+a_{j+1}}}{3^{j+1}}-\frac{2^{a_1+a_2+\ldots+a_j+2}}{3^{j+1}}+n
\end{equation*}
\begin{equation*}
    m =\frac{2^{a_1+a_2+\ldots+a_j}}{3^{j+1}}\times(2^{a_{j+1}}-2^2)+n
\end{equation*}
which is the equation we have used for the theorem. Now we have to show that $a_{j+1} = 2b3^j + 2, b \in \mathbb{N}_0$. By definition, both $n$ and $m$ are integers therefore the first term on the right-hand side of equation \ref{EQ: Additive Formula} must be an integer. Since $3^{j+1} \nmid 2^{a_1+a_2+\ldots+a_j}$ it is required that $3^{j+1} \mid 2^{a_{j+1}}-2^2$. This gives the congruence
\begin{equation*}
    2^{a_{j+1}}-2^2 \equiv 0 \pmod{3^{j+1}}
\end{equation*}
Plugging in $a_{j+1} = 2b3^j + 2$ yields
\begin{equation*}
    2^{2b3^j + 2}-2^2 \equiv 0 \pmod{3^{j+1}}
\end{equation*}
We can parenthesize with $2^2$ which gives
\begin{equation*}
    2^2(2^{2b3^j}-1) \equiv 0 \pmod{3^{j+1}}
\end{equation*}
Since $3^{j+1} \nmid 2^2$ our congruence is now
\begin{equation}
\label{EQ: Additive Formula Exponent Congruence}
    2^{2b3^j}-1 \equiv 0 \pmod{3^{j+1}}
\end{equation}
We will use proof by induction to prove this congruence. First, we look at $b=0$ which gives $2^{0}-1 \equiv 0 \pmod{3^{j+1}}$ which is correct. We also look at $b=1$ which gives 
\begin{equation}
\label{EQ: b=1 hali}
    2^{2\times3^j}-1 \equiv 0 \pmod{3^{j+1}}
\end{equation}
We will come to this in a moment. Next, we assume that $b=k$ is true therefore 
\begin{equation*}
    2^{2\times k\times3^j}-1 \equiv 0 \pmod{3^{j+1}}
\end{equation*}
Finally we look at $b=k+1$ which gives
\begin{equation*}
    2^{2\times(k+1)\times3^j}-1 \equiv 0 \pmod{3^{j+1}}
\end{equation*}
We can see that the expression $2^{2\times(k+1)\times3^j}-1 = 2^{2\times k\times3^j}\times2^{2\times1\times3^j}-1$. Adding $2^{2\times k\times3^j}-2^{2\times k\times3^j}$ to this gives
\begin{equation*}
    2^{2\times k\times3^j}\times2^{2\times 1\times3^j}-1+2^{2\times k\times3^j}-2^{2\times k\times3^j} \equiv 0 \pmod{3^{j+1}}
\end{equation*}
Parenthesizing with $2^{2\times k\times3^j}$ gives
\begin{equation*}
    2^{2\times k\times3^j}(2^{2\times 1\times 3^j}-1)+2^{2\times k\times 3^j}-1 \equiv 0 \pmod{3^{j+1}}
\end{equation*}
Here the second and third terms together is the same expression we had for $b=k$. For the first term, we have $2^{2\times 1\times 3^j}-1$ in parenthesis, which is the same expression we had for $b=1$ (equation \eqref{EQ: b=1 hali}). This tells us that we can prove the congruence \eqref{EQ: Additive Formula Exponent Congruence} by proving it is correct for $b=1$. In other words, we will prove by induction that the congruence
\begin{equation}
    2^{2\times3^j}-1 \equiv 0 \pmod{3^{j+1}}
\end{equation}
holds for every $j$. Starting with $j=1$ we have 
\begin{equation*}
    2^{2\times3^{1}}-1 \equiv 0 \pmod{3^{1+1}}
\end{equation*}
which gives $63 \equiv 0 \pmod{9}$ and is correct. Next, we look at $j=k$
\begin{equation*}
\label{EQ: j=k stepi}
    2^{2\times3^k}-1 \equiv 0 \pmod{3^{k+1}}
\end{equation*}
Assuming this is correct we look at $j=k+1$
\begin{equation*}
    2^{2\times3^{k+1}}-1 \equiv 0 \pmod{3^{k+1+1}}
\end{equation*}
We can write $2^{2\times3^{k+1}}$ as $2^{2\times3^{k}}\times2^{4\times3^{k}}$. Adding $2^{2\times3^k}-2^{2\times3^k}$ to this gives 
\begin{equation*}
    2^{2\times3^k}\times2^{4\times3^k}-1+2^{2\times3^k}-2^{2\times3^k} \equiv 0 \pmod{3^{k+1+1}}
\end{equation*}
Parenthesizing with $2^{2\times3^k}$ gives
\begin{equation*}
2^{2\times3^k}\times(2^{4\times3^k}-1)+2^{2\times3^k}-1 \equiv 0 \pmod{3^{k+1+1}}
\end{equation*}
Note that $2^{4\times3^k} = 2^{2\times3^k}\times 2^{2\times3^k}$ 
\begin{equation}
\label{EQ: Bir suru kuvvet var}
    2^{2\times3^k}\times(2^{2\times3^k}\times 2^{2\times3^k}-1)+2^{2\times3^k}-1 \equiv 0 \pmod{3^{k+1+1}}
\end{equation}
To show the next step clearly, we will denote $2^{2\times3^k} - 1$ as $A$, so $A = 2^{2\times3^k} - 1$. The expression at congruence \ref{EQ: Bir suru kuvvet var} then becomes
\begin{flalign*}
(A+1)((A+1)^2-1)+A &\equiv 0 \pmod{3^{k+1+1}}\\
(A+1)(A^2+2A+1-1)+A &\equiv 0 \pmod{3^{k+1+1}}\\
(A+1)(A^2+2A)+A &\equiv 0 \pmod{3^{k+1+1}}\\
A^3+2A^2+A^2+2A+A &\equiv 0 \pmod{3^{k+1+1}}\\
A^3+3A^2+3A &\equiv 0 \pmod{3^{k+1+1}}
\end{flalign*}
Opening $A$ back up and writing the congruence shows us that
\begin{equation}
    (2^{2\times3^k} - 1)^3+3(2^{2\times3^k} - 1)^2+3(2^{2\times3^k} - 1) \equiv 0 \pmod{3^{k+1+1}}
\end{equation}
Now from congruence \eqref{EQ: j=k stepi} we have $2^{2\times3^k}-1 \equiv 0 \pmod{3^{k+1}}$. From this we can definitely say that the congruence for the first term $(2^{2\times3^k} - 1)^3 \equiv 0 \pmod{3^{k+1+1}}$ holds, similarly the congruence for the second term $3(2^{2\times3^k} - 1)^2 \equiv 0 \pmod{3^{k+1+1}}$ holds and the last congruence $3(2^{2\times3^k} - 1) \equiv 0 \pmod{3^{k+1+1}}$ holds. Therefore it is proven that the congruence  \eqref{EQ: Additive Formula Exponent Congruence} holds for all $b \in \mathbb{N}_0$ and $j \in \mathbb{N}$.
\end{proof}
Sadly, this formula does not give every odd number $m$ with reduced Collatz sequence length $j+1$ when we have an odd number $n$ with reduced Collatz sequence length $j$. To demonstrate, we can look at the formula for numbers with reduced Collatz sequence length $1$ which we have shown in theorem \ref{THEO: Formula for Length 1 theorem}. Let $n \in \mathbb{O}_+$ with reduced Collatz sequence length 1 and $m \in \mathbb{O}_+$ with reduced Collatz sequence length $2$. Theorem \ref{THEO: Additive Formula} states that
\begin{equation*}
    m = \frac{2^{a_1}}{3^2}\times(2^{a_2}-2^2) + n
\end{equation*}
Theorem \ref{THEO: Formula for Length 1 theorem} stated that 
\begin{equation*}
    n = \frac{2^{2k+2}}{3^1} - \frac{2^0}{3^1}
\end{equation*}
So we have
\begin{equation*}
    m = \frac{2^{2k+2}}{3^2}\times(2^{a_2}-2^2) + \frac{2^{2k+2}}{3^1} - \frac{2^0}{3^1}
\end{equation*}
Theorem \ref{EQ: Additive Formula} states that $a_2 = 2b3^1+2$ so we have 
\begin{equation*}
    m = \frac{2^{2k+2}}{3^2}\times(2^{6b+2}-2^2) + \frac{2^{2k+2}}{3^1} - \frac{2^0}{3^1}
\end{equation*}
\begin{equation*}
    m = \frac{2^{2k+6b+4}}{3^2} - \frac{2^{2k+4}}{3^2} + \frac{2^{2k+2}}{3^1} - \frac{2^0}{3^1}
\end{equation*}
\begin{equation}
\label{EQ: Additive ile elde edilen length 2}
    m = \frac{2^{2k+6b+4}}{3^2} - \frac{2^{2k+2}}{3^2} - \frac{2^0}{3^1}
\end{equation}
However, theorem \ref{THEO: Formula for Length 2 theorem} states that
\begin{equation}
\label{EQ: Normal elde edilen length 2 ama asli}
    m = \frac{2^{a_1+6b+2(-1)^{a_1 \bmod{2}}}}{3^2}-\frac{2^{a_1}}{3^2}-\frac{2^0}{3^1}
\end{equation}
Now that we know $a_1 = 2k+2$ we get
\begin{equation}
\label{EQ: Normal elde edilen length 2}
    m = \frac{2^{2k+6b+4}}{3^2}-\frac{2^{2k+2}}{3^2}-\frac{2^0}{3^1}
\end{equation}
Equations \eqref{EQ: Additive ile elde edilen length 2} and \eqref{EQ: Normal elde edilen length 2} are same. The problem is, we have shown in theorem \ref{EQ: Formula for Length 2} that the equation \eqref{EQ: Normal elde edilen length 2 ama asli} is also correct. Equation \eqref{EQ: Normal elde edilen length 2 ama asli} yields more numbers with reduced Collatz sequence length 2. Why does this happen? It's because equation \eqref{EQ: Additive Formula}, when used with a number $n$ with reduced Collatz sequence length $j$ to obtain a number $m$ with reduced Collatz sequence length $j+1$, is constrained to the fact that $n$ is an odd integer. In the case of reduced Collatz sequence length 2, we were constrained to the fact that $n$ with reduced Collatz sequence length 1 was also valid. Having proved that the exponent must be an even number for reduced Collatz sequence length 1, we constrain ourselves to use even values for $a_1$ when we use \eqref{EQ: Additive Formula}. Theorem \ref{EQ: Formula for Length 2} clearly shows that it can be odd too. 

\subsection{Trivial Cycle}
\label{SECT: Trivial Cycle}
One pretty result of equation \eqref{EQ: Additive Formula} is when $a_{j+1} = 2$.
\begin{equation}
    m = \frac{2^{a_1+a_2+\ldots+a_j}}{3^{j+1}}\times(2^{2}-2^2)+n
\end{equation}
This yields $m = 0 + n$. The trajectory of $m$ is then $m\xrightarrow{a_1}C(m)\xrightarrow{a_2}C^2(m)\ldots\xrightarrow{a_{j-1}}C^{j-1}(m)\xrightarrow{a_j}C^j(m)\xrightarrow{2}1$. Looking at a trajectory such as $k \xrightarrow{2} 1$ we can see that $(3k+1)/2^2 = 1$ so $k=1$. We also know $m=n$. If $a_{j+1} = 2$ it just means that the trajectory ended at $1$ already and is looping. 

An interesting result of using the trivial cycle is that given a number $n$ with reduced Collatz sequence length $j$ we can find numbers with greater sequence lengths, not just $j+1$ but any length $j+k, k \in \mathbb{N}$.
\begin{corollary}
\label{COROL: Additive Formula using Trivial Cycle}
Let $n, m, k \in \mathbb{N}$ where $n$ has reduced Collatz sequence length $j$ and $m$ has reduced Collatz sequence length $j+k$.
\begin{equation}
\label{EQ: Additive Formula using Trivial Cycle}
    m = \frac{2^{a_1+a_2+\ldots+a_j+2(k-1)}}{3^{j+k}}\times(2^{a_{j+k}}-2^2)+n
\end{equation}
where $a_{j+k} = 2b3^{j+k-1} + 2, b \in \mathbb{N}_0$ and $a_{j+1} = a_{j+2} = \ldots = a_{j+k-1} = 2$ so for $k-1$ iterations we set the new exponent $2$ and therefore obtain a new formula using the same base $n$. The reduced Collatz sequences can be shown as
$$
n\xrightarrow{a_1}C(n)\xrightarrow{a_2}\ldots\xrightarrow{a_j}1
$$
$$
m \xrightarrow{a_1} C(m) \xrightarrow{a_2} \ldots \xrightarrow{a_j} C^j(m) 
\underbrace{\xrightarrow{2} C^{j+1}(m) \xrightarrow{2} \ldots \xrightarrow{2} C^{j+k-1}(m) }_{k-1\text{ times}}
\xrightarrow{a_{j+k}} 1
$$
\end{corollary}
Note how plugging $k=1$ in equation \eqref{EQ: Additive Formula using Trivial Cycle} basically yields equation \eqref{EQ: Additive Formula}. 
\begin{example}
Let us consider the number 5. We will demonstrate both equations \eqref{EQ: Additive Formula using Trivial Cycle} and \eqref{EQ: Additive Formula} with 5 as our base number ($n=5$). First let us look at the reduced Collatz sequence of 5 which is $5 \xrightarrow{4} 1$. It has reduced Collatz sequence length $1$. Equation \eqref{EQ: Additive Formula} becomes
\begin{equation}
    m = \frac{2^{4}}{3^{2}}\times(2^{6b + 2}-2^2)+5
\end{equation}
If $b=0$ we get 5 with the sequence $5 \xrightarrow{4} 1 \xrightarrow{2} 1$, which has looped over $1$ for one iteration and the new exponent is $2$ therefore
$$
5 = \frac{2^6}{3^2} - \frac{2^4}{3^2} - \frac{2^0}{3^1}
$$
If $b=1$ we get 
\begin{equation*}
m = \frac{2^{4}}{3^{2}}\times(2^{8}-2^2)+5
\end{equation*}
We find $m=453$ and the theory states that this number should have reduced Collatz sequence length $2$. Indeed it has, with the sequence $453 \xrightarrow{4} 85 \xrightarrow{8} 1$. 
$$
453 = \frac{2^{12}}{3^2} - \frac{2^4}{3^2} - \frac{2^0}{3^1}
$$
Finally we will show how we can obtain a number with reduced Collatz sequence length $3$ by using 5, which has reduced Collatz sequence length $1$. Corollary \ref{COROL: Additive Formula using Trivial Cycle} stated that we can obtain length $j+k$ with $k-1$ iterations with exponent $2$. We will do exactly that. In this example $j=1$ and we want $j+k=3$ so $k=2$. Plugging these in equation \eqref{EQ: Additive Formula using Trivial Cycle} we get
\begin{equation*}
m = \frac{2^{4+2(2-1)}}{3^{1+2}}\times(2^{2b3^{1+2-1} + 2}-2^2)+5    
\end{equation*}
\begin{equation*}
m = \frac{2^{6}}{3^{3}}\times(2^{18b + 2}-2^2)+5    
\end{equation*}
We will just give $b=1$ which results in
\begin{equation*}
m = \frac{2^{6}}{3^{3}}\times(2^{20}-2^2)+5   
\end{equation*}
This yields $m = 2485509$. It's reduced Collatz sequence is $2485509 \xrightarrow{4} 466033 \xrightarrow{2} 349525 \xrightarrow{20} 1$ and the fractional sum notation is
$$
2485509 = \frac{2^{26}}{3^3} - \frac{2^6}{3^3} - \frac{2^4}{3^2} - \frac{2^0}{3^1}
$$
\end{example}


%% file: Sections/3-IFSN.tex
\section{Intermediate Fractional Sum Notation}
So far we have only used fractional sum notation for Collatz trajectories that end with $1$. In this section we will introduce the Intermediate Fractional Sum Notation which does not have that requirement. Any Collatz trajectory is applicable.
\begin{definition}
\label{DEF: IFSN Definition}
Let $n, m \in \mathbb{O}_+$, $C^j(n) = m$. The reduced Collatz trajectory is
$$
n \xrightarrow{a_1} C^1(n) \xrightarrow{a_2}\ldots\xrightarrow{a_{j-1}}C^{j-1}(n) \xrightarrow{a_j} m
$$
Then we can write $n$ as
\begin{equation}
\label{DEF: IFSN}
    n = m\times \frac{2^{a_1+a_2+\ldots+a_j}}{3^j}-\frac{2^{a_1+a_2+\ldots+a_{j-1}}}{3^j}-\frac{2^{a_1+a_2+\ldots+a_{j-2}}}{3^{j-1}}-\ldots-\frac{2^{a_1}}{3^2}-\frac{2^{0}}{3^1}
\end{equation}
\end{definition}
Essentially, this notation and the fractional sum notation are same. Giving $m=1$ turns this into a fractional sum notation. 
\begin{example}
The reduced Collatz sequence of $7$ is
$$
7 \xrightarrow{1} 11 \xrightarrow{1} 17 \xrightarrow{2} 13 \xrightarrow{3} 5 \xrightarrow{4} 1
$$
So, the Collatz trajectory of $7$ up to $13$ is
$$
7 \xrightarrow{1} 11 \xrightarrow{1} 17 \xrightarrow{2} 13
$$
Writing Intermediate Fractional Sum Notation for $n=7$ and $m=C^3(n)=13$ gives us
$$
7=13\times\frac{2^4}{3^3}-\frac{2^2}{3^3}-\frac{2^1}{3^2}-\frac{2^0}{3^1}
$$
\end{example}
If there exists a number that does not reach 1 at some iteration, then it is either looping or diverging. So one may find that there is a number with a looping trajectory or a diverging trajectory to prove that this conjecture is false. The intermediate fractional sum notation is a tool we can use to study such trajectories.

%% file: Sections/3.1-Loops.tex
\subsection{Studying looping trajectories with intermediate fractional sum notation}
\label{SECT: Loops}
Suppose there is a number $n \in \mathbb{O}_+$ such that $C^j(n)=n$. The trajectory can be shown as
$$
n \xrightarrow{a_1} C(n) \xrightarrow{a_2} \ldots \xrightarrow{a_{j-1}} C^{j-1}(n) \xrightarrow{a_j} n
$$
This number is in a looping trajectory of length $j$. In fact, every number in this loop is looping (as the words suggest). The intermediate fractional sum notation for $n$ is
\begin{equation}
\label{EQ: IFSN for looping n}
    n = n\times \frac{2^{a_1+a_2+\ldots+a_j}}{3^j}-\frac{2^{a_1+a_2+\ldots+a_{j-1}}}{3^j}-\frac{2^{a_1+a_2+\ldots+a_{j-2}}}{3^{j-1}}-\ldots-\frac{2^{a_1}}{3^2}-\frac{2^{0}}{3^1}
\end{equation}
\begin{theorem}
\label{THEO: Looping n form}
Let $n \in \mathbb{O}_+$ and $C^j(n)=n$ with the trajectory $n \xrightarrow{a_1} C(n) \xrightarrow{a_2} \ldots \xrightarrow{a_{n-1}} C^{j-1}(n) \xrightarrow{a_j} n$. If there exists a number $m \in \mathbb{O}_+$ such that $C^j(m) = 1$ and the trajectory is $m \xrightarrow{a_1} C(m) \xrightarrow{a_2} \ldots \xrightarrow{a_{j-1}} C^{j-1}(m) \xrightarrow{a_j} 1$ then $n$ is in the form $3^j\times 2k + 1, k \in \mathbb{N}_0$.
\end{theorem}
\begin{proof}
We start by writing the intermediate fractional sum notation for $n$ which is
\begin{equation}
    n = n\times \frac{2^{a_1+a_2+\ldots+a_j}}{3^j}-\frac{2^{a_1+a_2+\ldots+a_{j-1}}}{3^j}-\frac{2^{a_1+a_2+\ldots+a_{j-2}}}{3^{j-1}}-\ldots-\frac{2^{a_1}}{3^2}-\frac{2^{0}}{3^1}
\end{equation}
In the first term on the right-hand side of the equation, we write $n+1-1$ instead of $n$, which gives us
\begin{equation*}
    n = (n+1-1) \times \frac{2^{a_1+a_2+\ldots+a_j}}{3^j}-\frac{2^{a_1+a_2+\ldots+a_{j-1}}}{3^j}-\frac{2^{a_1+a_2+\ldots+a_{j-2}}}{3^{j-1}}-\ldots-\frac{2^{a_1}}{3^2}-\frac{2^{0}}{3^1}
\end{equation*}
\begin{equation}
\label{EQ: looplanmis n+1-1}
    n = (n-1) \times \frac{2^{a_1+\ldots+a_j}}{3^j}+\frac{2^{a_1+\ldots+a_j}}{3^j}-\frac{2^{a_1+\ldots+a_{j-1}}}{3^j}-\frac{2^{a_1+\ldots+a_{j-2}}}{3^{j-1}}-\ldots-\frac{2^{a_1}}{3^2}-\frac{2^{0}}{3^1}
\end{equation}
If there is a number $m \in \mathbb{O}_+$ with the trajectory $m \xrightarrow{a_1} C(m) \xrightarrow{a_2} \ldots \xrightarrow{a_{j-1}} C^{j-1}(m) \xrightarrow{a_j} 1$ then the fractional sum notation for $m$ is
\begin{equation*}
    m = \frac{2^{a_1+a_2+\ldots+a_j}}{3^j}-\frac{2^{a_1+a_2+\ldots+a_{j-1}}}{3^j}-\frac{2^{a_1+a_2+\ldots+a_{j-2}}}{3^{j-1}}-\ldots-\frac{2^{a_1}}{3^2}-\frac{2^{0}}{3^1}
\end{equation*}
Then equation \eqref{EQ: looplanmis n+1-1} becomes
\begin{equation}
\label{EQ: looplanmis n ve m var}
    n = (n-1) \times \frac{2^{a_1+a_2+\ldots+a_j}}{3^j} + m
\end{equation}
We define $n$ and $m$ to be odd numbers, therefore the first term on the right-hand side must be an integer. Since $3^j \nmid 2^{a_1+a_2+\ldots+a_j}$ we must have $3^j \mid n-1$. This gives the congruence
\begin{equation}
    n - 1 \equiv 0 \pmod{3^j}
\end{equation}
Looking at this congruence we can say that $n = 3^j \times k+1, k \in \mathbb{N}_0$. There is one more thing to consider, $n$ is an odd number. For odd values of $k$ the expression $3^j \times k+1$ yields even numbers because $(\text{odd}\times\text{odd}+\text{odd})$ results in an $\text{even}$ number. Since $n$ must be odd we can not have odd values for $k$. For even values of $k$ the expression $3^j \times k+1$ yields odd numbers. Rewriting the expression with $2k$ instead of $k$ gives $n = 3^j\times 2k + 1, k \in \mathbb{N}_0$.
\end{proof}
\begin{corollary}
Under the same conditions mentioned in theorem \ref{THEO: Looping n form}, since $n = 3^j\times 2k + 1$ we can say $n \equiv 1 \pmod{6}$. Modulo 6 turns out to be quite important, as we will see in section \ref{SECT: Reverse Formulation}.
\end{corollary}
So far we only have one known loop: $1 \xrightarrow{2} 1$. Plugging $n=1$ in equation \eqref{EQ: looplanmis n ve m var} we get
\begin{equation*}
    1 = (1-1)\times\frac{2^2}{3^1} + m
\end{equation*}
which shows that $m=1$, which is consistent with the condition stated in theorem \ref{THEO: Looping n form}. In fact, a looping trajectory
$$
1\underbrace{ \xrightarrow{2} 1 \xrightarrow{2} \ldots \xrightarrow{2} 1 \xrightarrow{2} 1}_{j \text{ iterations}}
$$
is consistent with the theorem. All exponents are $2$ and if there are $j$ iterations equation \eqref{EQ: looplanmis n ve m var} becomes
\begin{equation*}
    1 = (1-1)\frac{2^{2j}}{3^j} + m
\end{equation*}
which still gives $m=1$.
\subsubsection{Loop with length 2}

It has been proved by J. L. Simons \cite{simons} that there are no loops with length 2 other than $1 \xrightarrow{2} 1 \xrightarrow{2} 1$, the proof used upper bounds and lower bounds, following R. P. Steiner's method \cite{steiner}. We can use intermediate fractional sum notation to quickly prove there are no loops with length 2 other than $1 \xrightarrow{2} 1 \xrightarrow{2} 1$ too. Consider the trajectory $n \xrightarrow{a_1} C(n) \xrightarrow{a_2} n$. The intermediate fractional sum notation is
\begin{equation}
    n = n\times \frac{2^{a_1+a_2}}{3^2} - \frac{2^{a_1}}{3^2} - \frac{2^0}{3^1}
\end{equation}
This gives us the equation $9n = n\times2^{a_1+a_2} - 2^{a_1} - 3$. Leaving $n$ alone,
\begin{equation}
\label{EQ: N alone in loop length 2}
    n = \frac{2^{a_1}+3}{2^{a_1+a_2}-9}
\end{equation}
This is equal to
\begin{equation*}
    n = \frac{2^{a_1}+3}{2^{a_1}(2^{a_2}+1-1)-12+3}
\end{equation*}
\begin{equation*}
    n = \frac{2^{a_1}+3}{2^{a_1}+3+2^{a_1}(2^{a_2}-1)-12}
\end{equation*}
Now remember that $n$ is an integer, therefore numerator is greater than or equal to denominator. This means that $2^{a_1}+3 \geq 2^{a_1}+3+2^{a_1}(2^{a_2}-1)-12$ which reduces to $ 12 \geq 2^{a_1+a_2} - 2^{a_1} $. Since $n$ is a positive odd integer, in equation \eqref{EQ: N alone in loop length 2} at the denominator we get the inequality $2^{a_1+a_2} > 9$. This tells us $a_1+a_2 > 3$. Also remember that $a_1, a_2 \in \mathbb{N}$. If $a_1=1$ we get $14 \geq 2^{1+a_2}$, which is possible for $a_2 < 3$, but then adding $a_1$ to this gives $1+a_2 < 4$ which contradicts the inequality $a_1+a_2 > 3$. If $a_1 = 2$ we get $16 \geq 2^{2+a_2}$, which holds for $a_2=1$ and $a_2=2$ but $2+1 < 3$ is false therefore only valid value is $a_2 = 2$. If $a_1 = 3$ we get $20 \geq 2^{3+a_2}$. Only valid value is $a_2=1$. If $a_1 = 4$ we get $28 \geq 2^{4+a_2}$. There are no valid values for $a_2$ after this point. 2 pairs of $(a_1, a_2) $ satisfy both inequalities: $(2, 2)$ and $(3, 1)$. Looking at pair $(3, 1)$ if we plug it in equation \eqref{EQ: N alone in loop length 2} we get
\begin{equation*}
    n = \frac{2^{3}+3}{2^{4}-9} = \frac{11}{7} 
\end{equation*}
which does not comply with $n \in \mathbb{O}_+$ as $11/7$ is not an integer. The other pair is $(2,2)$ which gives
\begin{equation*}
    n = \frac{2^{2}+3}{2^{4}-9} = \frac{7}{7} 
\end{equation*}
therefore $n=1$. This shows that the only loop of length 2 is $1 \xrightarrow{2} 1\xrightarrow{2} 1$.
\subsubsection{Loop with greater lengths}
\label{SECT: Loop with  greater length}
Can we do the same thing for a loop with length $3$? The trajectory would be $n \xrightarrow{a_1} C(n) \xrightarrow{a_2} C^2(n) \xrightarrow{a_3} n$.
\begin{equation}
    n = n\times \frac{2^{a_1+a_2+a_3}}{3^3} - \frac{2^{a_1+a_2}}{3^3} - \frac{2^{a_1}}{3^2} - \frac{2^0}{3^1}
\end{equation}
Leaving $n$ alone gives
\begin{equation}
    n = \frac{2^{a_1+a_2} + 3\times 2^{a_1} + 9}{2^{a_1+a_2+a_3} - 27}
\end{equation}
This is equal to 
\begin{equation*}
    n = \frac{2^{a_1+a_2} + 3\times 2^{a_1} + 9}{2^{a_1+a_2} + 3\times 2^{a_1} + 9 + 2^{a_1+a_2}(2^{a_3}-1-\frac{3}{2^{a_2}})-36}
\end{equation*}
We can try the inequality trick on $2^{a_1+a_2} + 3\times 2^{a_1} + 9 \geq 2^{a_1+a_2} + 3\times 2^{a_1} + 9 + 2^{a_1+a_2}(2^{a_3}-1-\frac{3}{2^{a_2}})-36$ together with $2^{a_1+a_2+a_3} - 27 > 0$, but it seems things get too complicated at this point. Considering $n \xrightarrow{a_1} C(n) \xrightarrow{a_2} \ldots \xrightarrow{a_j} n$ gives equation \eqref{EQ: IFSN for looping n}. Applying the same procedure to leave $n$ alone yields
\begin{equation}
    n = \frac{\displaystyle \sum_{i=1}^{j}\left(3^{i-1}\times2^{\sum_{k=1}^{j-i}a_k}\right)}{2^{\sum_{i=1}^{j}a_i}-3^j}
\end{equation}
If there is an odd integer $n$ with loop length $j$ then it must satisfy this equation. Not the most malleable formula, but it is there!

%% file: Sections/3.2-Divergence.tex
\subsection{Studying diverging trajectories with intermediate fractional sum notation}
F.C. Motta et al. \cite{fcmotta} worked on diverging trajectories in their paper, they had a formula that gives an upper bound for a monotonically increasing trajectory. A monotonically increasing trajectory is a trajectory where on each iteration the result gets bigger, for example
$$
n \xrightarrow{1} C(n) \xrightarrow{1} C^2(n) \xrightarrow{1} \ldots
$$
We can use the intermediate fractional sum notation to find numbers which monotonically increase for $j$ iterations, with the trajectory
$$
n \xrightarrow{1} C(n) \xrightarrow{1} \ldots \xrightarrow{1} C^j(n)
$$
The intermediate fractional sum notation is
\begin{equation}
    n = C^j(n)\times\frac{2^{j}}{3^j}-\frac{2^{j-1}}{3^j}-\frac{2^{j-2}}{3^{j-1}}-\ldots-\frac{2^{1}}{3^2}-\frac{2^{0}}{3^1}
\end{equation}
\begin{equation*}
    n = C^j(n)\times\frac{2^{j}}{3^j}-\frac{1}{3}\times\left(\frac{2^{j-1}}{3^{j-1}}+\frac{2^{j-2}}{3^{j-2}}+\ldots+\frac{2^{1}}{3^1}+1\right)
\end{equation*}
Inside the parentheses we have a finite geometric series
\footnote{The sum of first $n$ terms in a geometric series is given by
$$
\frac{a_1\times(1-r^n)}{1-r}
$$
where $a_1$ is the first term and $r$ is the common ratio, $-1<r<1$. In our geometric series $a_1=1$ and $r=\frac{2}{3}$}%
.Writing the sum of $j$ terms for that geometric series gives us
\begin{equation*}
n = C^j(n) \times\frac{2^{j}}{3^j}-\frac{1}{3}\times\left(\frac{1-\frac{2^j}{3^j}}{1-\frac{2}{3}}\right)    
\end{equation*}
Leaving $C^j(n)$ alone we get
\begin{equation}
\label{EQ: Divergent Trajectory}
    C^j(n) = \left(\frac{3}{2}\right)^j(n+1)-1
\end{equation}
Other than a minor symbol difference, this is the same equation seen in the paper of F.C. Motta et al. (\cite{fcmotta} equation (1)).
\begin{example}
One can construct monotonically increasing trajectories using equation \eqref{EQ: Divergent Trajectory}. Giving $n=2^{j}\times k-1, k\in\mathbb{N}$ yields a monotonically increasing trajectory with $j$ odd numbers. $7=2^3-1$, monotonically increasing trajectory with $3$ numbers:$7 \xrightarrow{1} 11 \xrightarrow{1} 17$. $15=2^4-1$, monotonically increasing trajectory with $4$ numbers: $15 \xrightarrow{1} 23 \xrightarrow{1} 35 \xrightarrow{1} 53$. $39=5\times2^3-1$, monotonically increasing trajectory with $3$ numbers: $39 \xrightarrow{1} 59 \xrightarrow{1} 89$. $127=2^7-1$, monotonically increasing trajectory with $7$ numbers:
$127 \xrightarrow{1} 191 \xrightarrow{1} 287 \xrightarrow{1} 431 \xrightarrow{1} 647\xrightarrow{1}  971 \xrightarrow{1} 1457$ and so on...
\end{example}

%% file: Sections/4-ReverseFormula.tex
\section{Reverse Formulation of the problem}
\label{SECT: Reverse Formulation}
A widely used reformulation of the problem is using the reverse reduced Collatz function $R : \mathbb{O}_+ \xrightarrow{} \mathbb{O}_+$ defined by
\begin{equation}
\label{EQ: Reverse Formula}
    R(n) = \frac{2^x\times n -1}{3}
\end{equation}
where $n \in \mathbb{O}_+$ and $x \in \mathbb{N}$. The problem so far was asking whether any odd number $n$ reach 1 at some iteration $C^i(n)$. We can ask the same by saying that for every odd number $n$ there is an iteration $R^i(1) = n$. Regarding $R$ function there are 3 lemmas, relating to modulo 6.
\begin{lemma}
\label{LEMMA: 6a+1 ise x even}
If $n \equiv 1 \pmod{6}$ then $x$ is an even number.
\end{lemma}
\begin{proof}
We can write $n = 6a+1, a \in \mathbb{N}_0$ and $x = 2b+2, b \in \mathbb{N}_0$. Equation \eqref{EQ: Reverse Formula} becomes
\begin{equation}
    R(6a+1) = \frac{(6a+1)\times 2^{2b+2} -1}{3}
\end{equation}
For this function to be valid the congruence $(6a+1)\times 2^{2b+2}  -1 \equiv 0 \pmod{3}$ must hold. 
\begin{equation*}
     6a\times 2^{2b+2} + 2^{2b+2} -1 \equiv 0 \pmod{3}
\end{equation*}
Clearly $3 \mid 2^{2b+2}\times6a$ and $3 \mid 2^{2b+2} -1$ was shown in lemma \ref{LEMMA: Length 1 Even power}.
\end{proof}

\begin{lemma}
\label{LEMMA: 6a+3 ise x yok}
If $n \equiv 3 \pmod{6}$ there are no possible values for $x$.  
\end{lemma}
\begin{proof}
We can write $n = 6a+3, a \in \mathbb{N}_0$. Equation \eqref{EQ: Reverse Formula} becomes
\begin{equation}
    R(6a+3) = \frac{(6a+3)\times 2^{x}  -1}{3}
\end{equation}
For this function to be valid the congruence $(6a+3)\times 2^{x}  -1 \equiv 0 \pmod{3}$ must hold. $3 \mid (6a+3)\times 2^{x} $ but $3 \nmid -1$ therefore the congruence does not hold for any $x$.
\end{proof}

\begin{lemma}
\label{LEMMA: 6a+5 ise x odd}
If $n \equiv 5 \pmod{6}$ then $x$ is an odd number.  
\end{lemma}
\begin{proof}
We can write $n = 6a+5, a \in \mathbb{N}_0$ and $x = 2b+1, b \in \mathbb{N}_0$. Equation \eqref{EQ: Reverse Formula} becomes
\begin{equation}
    R(6a+5) = \frac{(6a+5)\times 2^{2b+1} -1}{3}
\end{equation}
For this function to be valid the congruence $(6a+5)\times2^{2b+1}  - 1 \equiv 0 \pmod{3}$ must hold. 
\begin{equation*}
     (6a+3)\times2^{2b+1} + 2\times2^{2b+1} -1 \equiv 0 \pmod{3}
\end{equation*}
Clearly $3 \mid (6a+3)\times2^{2b+1}$ and $3 \mid 2\times2^{2b+1} -1$ was shown in lemma \ref{LEMMA: Length 1 Even power}.
\end{proof}
These 3 lemmas show us $R(n)$ can only be used when $n \equiv 1 \pmod{6}$ or $n \equiv 5 \pmod{6}$. 
\begin{definition}
The function $R_1 : \mathbb{N}_0\times\mathbb{N}_0 \xrightarrow{} \mathbb{O}_{R1}$ is defined as
\begin{equation}
    R_1(a,b) = \frac{(6a+1)\times 2^{2b+2}-1}{3}
\end{equation}
$\mathbb{O}_{R1} \subset \mathbb{O}_+$. This function acts as a reverse reduced Collatz iteration on a number $n = 6a+1$. 
\end{definition}
\begin{definition}
The function $R_5 : \mathbb{N}_0\times\mathbb{N}_0 \xrightarrow{} \mathbb{O}_{R5}$ is defined as
\begin{equation}
    R_5(a,b) = \frac{(6a+5)\times 2^{2b+1}-1}{3}
\end{equation}
$\mathbb{O}_{R5} \subset \mathbb{O}_+$. This function acts as a reverse reduced Collatz iteration on a number $n = 6a+5$. 
\end{definition}
We will soon show that $\mathbb{O}_{R1} \cap \mathbb{O}_{R5} = \varnothing$ and $\mathbb{O}_{R1} \cup \mathbb{O}_{R5} = \mathbb{O}_+$ but to do that we need to introduce one more function and prove several lemmas. 
\begin{definition}
The function $X : \mathbb{N}_0 \times \mathbb{N}_0 \xrightarrow{} \mathbb{O}_X$ is defined as
\begin{equation}
    X(a,b) = 2^{2b+3}a + \frac{2^{2b+4}-1}{3}
\end{equation}
$\mathbb{O}_X \subset \mathbb{O}_+$. We will soon see why this function is important.
\end{definition}
\begin{definition}
We need to provide one more useful notation for our proofs. The set of all possible values function $R_1(a,b)$ can produce where $a \in \mathbb{N}_0$ and $b \in \{0, 1, \ldots, k\}$ will be shown as $\mathbb{O}_{R1, k}$. The set of all possible values function $R_5(a,b)$ can produce where $a \in \mathbb{N}_0$ and $b \in \{0, 1, \ldots, k\}$ will be shown as $\mathbb{O}_{R5, k}$. The set of all possible values function $X(a,b)$ can produce where $a \in \mathbb{N}_0$ and $b \in \{0, 1, \ldots, k\}$ will be shown as $\mathbb{O}_{X, k}$. If $k$ is infinity, it basically means $b \in \mathbb{N}_0$. In that case we do not write infinity explicitly, for example $\mathbb{O}_{R1, \infty}$ is wrong, we use $\mathbb{O}_{R1}$ as defined before.
\end{definition} 

\begin{lemma}
\label{LEMMA: ODD = R1 + R5 + X}
$\mathbb{O}_{R1, 0} \cup \mathbb{O}_{R5, 0} \cup \mathbb{O}_{X, 0} = \mathbb{O}_+$ and $\mathbb{O}_{R1, 0}, \mathbb{O}_{R5, 0}, \mathbb{O}_{X, 0}$ are disjoint sets.
\end{lemma}
\begin{proof}
Let us look at functions $R_1(a,0), R_5(a,0)$ and $X(a,0)$
\begin{equation*}
    R_1(a,0) = \frac{(6a+1)\times 2^{2}-1}{3} = 8a+1
\end{equation*}
\begin{equation*}
    R_5(a,0) = \frac{(6a+5)\times 2^{1}-1}{3} = 4a+3
\end{equation*}
\begin{equation*}
    X(a,0) = 2^{3}a + \frac{2^{4}-1}{3} = 8a + 5
\end{equation*}
Let $\bar{r}$ denote the set of positive odd integers congruent to $r$ modulo $8$, $\bar{r} = \{n \in \mathbb{O}_+ \mid n \equiv r \pmod{8}\}$. It is evident that $\mathbb{O}_{R1, 0} = \bar{1}$, $ \mathbb{O}_{R5, 0} = \bar{3} \cup \bar{7}$ and $\mathbb{O}_{X, 0} = \bar{5}$. We are able to say that $\bar{1} \cup \bar{3} \cup \bar{5} \cup \bar{7} = \mathbb{O}_+$ as well as see that $\bar{1}, \bar{3} , \bar{5} , \bar{7}$ are disjoint.
\end{proof}
\begin{lemma}
\label{LEMMA: X = R1 + R5 + Xnew}
$(\mathbb{O}_{X, k} \setminus \mathbb{O}_{X, k-1}) = (\mathbb{O}_{R1, k+1} \setminus \mathbb{O}_{R1, k}) \cup (\mathbb{O}_{R5, k+1} \setminus \mathbb{O}_{R5, k}) \cup (\mathbb{O}_{X, k+1} \setminus \mathbb{O}_{X, k})$
\end{lemma}
\begin{proof}
Now this may appear to be a rather complicated lemma, but it is not, the notation is. We can elucidate it a bit: Consider an arbitrary positive integer $k$ and $a \in \mathbb{N}_0$, the lemma states that the set of all possible values produced by $X(a,k)$ is equal to the union of sets all possible values produced by $R_1(a,k+1)$, $R_5(a,k+1)$ and $X(a,k+1)$. We just used the set difference to achieve it. To show this, consider the cases of $a$ in modulo 4: it can be 0, 1, 2 or 3. We can say that the set of all possible values produced by $X(a,k)$ is then equal to the union of sets all possible values produced by $X(4a,k)$, $X(2a+1,k)$ and $X(4a+2,k)$. Let us look at the equation $X(4a,k) = R_1(a,k+1)$.
\begin{flalign*}
        2^{2k+3}(4a)+\frac{2^{2k+4}-1}{3} &= \frac{(6a+1)2^{2(k+1)+2}-1}{3} \\
        3\times2^{2k+3}(2^2a)+2^{2k+4}-1&=(6a+1)2^{2k+4}-1 \\
        6a2^{2k+4}+2^{2k+4}&=(6a+1)2^{2k+4} \\
        (6a+1)2^{2k+4}&=(6a+1)2^{2k+4} \\
        1&=1
\end{flalign*}
It is indeed correct, next let us look at the equation $X(2a+1,k) = R_5(a,k)$.
\begin{flalign*}
       2^{2k+3}(2a+1)+\frac{2^{2k+4}-1}{3}&=\frac{(6a+5)2^{2(k+1)+1}-1}{3} \\
    3\times2^{2k+3}(2a+1)+2^{2x+4}-1&=(6a+5)2^{2(k+1)+1}-1 \\
    6a2^{2k+3}+3\times2^{2k+3}+2^{2k+4}&=6a2^{2k+3}+5\times2^{2k+3} \\
    3\times2^{2k+3}+2\times2^{2k+3}&=5\times2^{2k+3} \\
    5\times2^{2k+3}&=5\times2^{2k+3} \\
    1&=1
\end{flalign*}
This is also correct. Finally, we look at the equation $X(4a+2, k) = X(a, k+1)$.
\begin{flalign*}
    2^{2k+3}(4a+2)+\frac{2^{2k+4}-1}{3}&=2^{2(k+1)+3}a+\frac{2^{2(k+1)+4}-1}{3} \\
    2^{2k+5}a+2^{2k+4}+\frac{2^{2k+4}-1}{3}&=2^{2k+5}a+\frac{2^{2k+2+4}-1}{3} \\
    3\times2^{2x+4}+2^{2x+4}-1&=2^{2x+6}-1 \\
    4\times2^{2x+4}&=4\times2^{2x+4} \\
    1&=1
\end{flalign*}
\end{proof}
\begin{lemma}
\label{LEMMA: Unifying sets}
$(\mathbb{O}_{R1,k+x} \setminus \mathbb{O}_{R1,k}) \cup (\mathbb{O}_{R1,k} \setminus \mathbb{O}_{R1,k-y}) = \mathbb{O}_{R1,k+x} \setminus \mathbb{O}_{R1,k-y}$. This is also true for $\mathbb{O}_{R5}$ and $\mathbb{O}_{X}$.
\end{lemma}
\begin{proof}
Consider $(\mathbb{O}_{R1,k+x} \setminus \mathbb{O}_{R1,k}) \cup (\mathbb{O}_{R1,k} \setminus \mathbb{O}_{R1,k-y})$ and $a \in \mathbb{N}_0$. $\mathbb{O}_{R1,k+x} \setminus \mathbb{O}_{R1,k}$ is the set of all possible results produced by $R_1(a,k+x), R_1(a,k+x-1), \ldots, R_1(a,k+1)$. Similarly $\mathbb{O}_{R1,k} \setminus \mathbb{O}_{R1,k-y}$ is the set of all possible results produced by $R_1(a,k), R_1(a,k-1), R_1(a,k-2), \ldots, R_1(a,k-y+1)$. As a result, $(\mathbb{O}_{R1,k+x} \setminus \mathbb{O}_{R1,k}) \cup (\mathbb{O}_{R1,k} \setminus \mathbb{O}_{R1,k-y})$ is the set of all possible results produced by $R_1(a,k+x), R_1(a,k+x-1), \ldots, R_1(a,k+1), R_1(a,k), R_1(a,k-1), R_1(a,k-2), \ldots, R_1(a,k-y+1)$, therefore we can see that
\begin{equation}
    (\mathbb{O}_{R1,k+x} \setminus \mathbb{O}_{R1,k}) \cup (\mathbb{O}_{R1,k} \setminus \mathbb{O}_{R1,k-y}) = \mathbb{O}_{R1,k+x} \setminus \mathbb{O}_{R1,k-y}
\end{equation}
This applies to the functions $R_5$ and $X$ too, by definition. Therefore
\begin{equation}
    (\mathbb{O}_{R5,k+x} \setminus \mathbb{O}_{R5,k}) \cup (\mathbb{O}_{R5,k} \setminus \mathbb{O}_{R5,k-y}) = \mathbb{O}_{R5,k+x} \setminus \mathbb{O}_{R5,k-y}
\end{equation}
and
\begin{equation}
    (\mathbb{O}_{X,k+x} \setminus \mathbb{O}_{X,k}) \cup (\mathbb{O}_{X,k} \setminus \mathbb{O}_{X,k-y}) = \mathbb{O}_{X,k+x} \setminus \mathbb{O}_{X,k-y}
\end{equation}
\end{proof}

\begin{theorem}
\label{THEO: R1 U R5 = O}
$\mathbb{O}_{R1} \cap \mathbb{O}_{R5} = \varnothing$ and $\mathbb{O}_{R1} \cup \mathbb{O}_{R5} = \mathbb{O}_+$
\end{theorem}
\begin{proof}
Lemma \ref{LEMMA: X = R1 + R5 + Xnew} shows us that when $k=0$ we have 
\begin{equation*}
    (\mathbb{O}_{X, 0} \setminus \mathbb{O}_{X, 0-1}) = (\mathbb{O}_{R1, 0+1} \setminus \mathbb{O}_{R1, 0}) \cup (\mathbb{O}_{R5, 0+1} \setminus \mathbb{O}_{R5, 0}) \cup (\mathbb{O}_{X, 0+1} \setminus \mathbb{O}_{X, 0})
\end{equation*}
$\mathbb{O}_{X, -1}$ is not possible since $b \geq 0$ so $\mathbb{O}_{X, -1} = \varnothing$. Therefore
\begin{equation}
\label{EQ: x=0 kume unionlari}
    \mathbb{O}_{X, 0} = (\mathbb{O}_{R1, 1} \setminus \mathbb{O}_{R1, 0}) \cup (\mathbb{O}_{R5, 1} \setminus \mathbb{O}_{R5, 0}) \cup (\mathbb{O}_{X, 1} \setminus \mathbb{O}_{X, 0})
\end{equation}
For $k=1$ we have
\begin{equation*}
    (\mathbb{O}_{X, 1} \setminus \mathbb{O}_{X, 0}) = (\mathbb{O}_{R1, 2} \setminus \mathbb{O}_{R1, 1}) \cup (\mathbb{O}_{R5, 2} \setminus \mathbb{O}_{R5, 1}) \cup (\mathbb{O}_{X, 2} \setminus \mathbb{O}_{X, 1})
\end{equation*}
Writing $(\mathbb{O}_{X, 1} \setminus \mathbb{O}_{X, 0})$ in equation \eqref{EQ: x=0 kume unionlari}
\begin{equation*}
    \mathbb{O}_{X, 0} = (\mathbb{O}_{R1, 1} \setminus \mathbb{O}_{R1, 0}) \cup (\mathbb{O}_{R5, 1} \setminus \mathbb{O}_{R5, 0}) \cup (\mathbb{O}_{R1, 2} \setminus \mathbb{O}_{R1, 1}) \cup (\mathbb{O}_{R5, 2} \setminus \mathbb{O}_{R5, 1}) \cup (\mathbb{O}_{X, 2} \setminus \mathbb{O}_{X, 1})
\end{equation*}
Continuing this process yields the expression below
\begin{flalign*}
    \mathbb{O}_{X, 0} = &(\mathbb{O}_{R1, 1} \setminus \mathbb{O}_{R1, 0}) \cup (\mathbb{O}_{R5, 1} \setminus \mathbb{O}_{R5, 0}) \\
    \cup & (\mathbb{O}_{R1, 2} \setminus \mathbb{O}_{R1, 1}) \cup (\mathbb{O}_{R5, 2} \setminus \mathbb{O}_{R5, 1}) \\
    \cup & (\mathbb{O}_{R1, 3} \setminus \mathbb{O}_{R1, 2}) \cup (\mathbb{O}_{R5, 3} \setminus \mathbb{O}_{R5, 2}) \\
    & \vdots \\
    \cup & (\mathbb{O}_{R1, k} \setminus \mathbb{O}_{R1, k-1}) \cup (\mathbb{O}_{R5, k} \setminus \mathbb{O}_{R5, k-1}) \\
    \cup & (\mathbb{O}_{X, k} \setminus \mathbb{O}_{X, k-1})
\end{flalign*}
Applying lemma \ref{LEMMA: Unifying sets} on this yields
\begin{equation*}
    \mathbb{O}_{X, 0} = (\mathbb{O}_{R1, k} \setminus \mathbb{O}_{R1, 0}) \cup (\mathbb{O}_{R5, k} \setminus \mathbb{O}_{R5, 0}) \cup (\mathbb{O}_{X, k} \setminus \mathbb{O}_{X, k-1})
\end{equation*}
We can infer that
\begin{equation*}
    \mathbb{O}_{X, 0} = (\mathbb{O}_{R1} \setminus \mathbb{O}_{R1, 0}) \cup (\mathbb{O}_{R5} \setminus \mathbb{O}_{R5, 0})
\end{equation*}
as $k$ goes to infinity. Lemma \ref{LEMMA: ODD = R1 + R5 + X} states that
\begin{equation*}
    \mathbb{O}_{R1, 0} \cup \mathbb{O}_{R5, 0} \cup \mathbb{O}_{X, 0} = \mathbb{O}_+
\end{equation*}
Now that we have shown $\mathbb{O}_{X, 0} = (\mathbb{O}_{R1} \setminus \mathbb{O}_{R1, 0}) \cup (\mathbb{O}_{R5} \setminus \mathbb{O}_{R5, 0})$ we can use it in lemma \ref{LEMMA: ODD = R1 + R5 + X} to have
\begin{equation*}
    \mathbb{O}_{R1, 0} \cup \mathbb{O}_{R5, 0} \cup (\mathbb{O}_{R1} \setminus \mathbb{O}_{R1, 0}) \cup (\mathbb{O}_{R5} \setminus \mathbb{O}_{R5, 0}) = \mathbb{O}_+
\end{equation*}
Since $(\mathbb{O}_{R1} \setminus \mathbb{O}_{R1, 0}) \cup \mathbb{O}_{R1, 0} = \mathbb{O}_{R1}$ and $(\mathbb{O}_{R5} \setminus \mathbb{O}_{R5, 0}) \cup \mathbb{O}_{R5, 0} = \mathbb{O}_{R5}$ we finally have
\begin{equation*}
    \mathbb{O}_{R1} \cup \mathbb{O}_{R5} = \mathbb{O}_+
\end{equation*}
\end{proof}
This proof is a safety assurance that we can use the functions $R_1$ and $R_5$ to study reverse reformulation of the problem. It is important to understand what these functions mean. $R_1(a,b)$ is a reverse iteration of an odd number $n$ such that $n \equiv 1 \pmod{6}$ where the number is multiplied by $2^{2b+2}$ and then 1 is subtracted from this, finally being divided by 3. Let the final result be $m$. See that
\begin{equation*}
    m = \frac{2^{2b+2}\times n - 1}{3}
\end{equation*}
Leaving $n$ alone,
\begin{equation*}
    \frac{3m+1}{2^{2b+2}} = n
\end{equation*}
So actually $C(m) = n$ and $m = R(n)$. Similarly for $R_5$ function we work with $n$ such that $n \equiv 5 \pmod{6}$ and the result is
\begin{equation*}
    m = \frac{2^{2b+1}\times n - 1}{3}
\end{equation*}
Leaving $n$ alone,
\begin{equation*}
    \frac{3m+1}{2^{2b+1}} = n
\end{equation*}
Again, $C(m) = n$ and $m = R(n)$. The proof tells us that every odd number is a result of reverse iteration of some odd number. This is a rather obvious statement, but a required one for the function to be valid. Another thing this proof shows is that every odd number is the result of a reverse iteration of exactly one odd number. Using two functions for cases of modulo 6 can also be seen in K. Andersson's work \cite{kerstin}.

%% file: Sections/4.1-Remainders.tex
\subsection{Determining modulo 6 after a single reverse iteration}
In this section we present a modified version of $R_1$ and $R_5$, namely $RR_1$ and $RR_5$. Before we continue we need to provide several lemmas.
\begin{lemma}
\label{LEMMA: 2^2x+1 equiv 2 mod 6}
$2^{2x+1} \equiv 2 \pmod{6}$ where $x \in \mathbb{N}_0$.
\end{lemma}
\begin{proof}
At $x=0$ we get $2 \equiv 2 \pmod{6}$, at $x=1$ we get $8\equiv 2\pmod{6}$. Assume that for $x=k$ the congruence 
\begin{equation*}
    2^{2k+1} \equiv 2 \pmod{6}    
\end{equation*}
is correct. Looking at $x = k + 1$ we have
\begin{equation*}
    2^{2(k+1)+1} \equiv 2 \pmod{6}
\end{equation*}
Adding $2^{2k+1}-2^{2k+1}$ to this expression yields
\begin{equation*}
    2^{2k+2+1} + 2^{2k+1} - 2^{2k+1} \equiv 2 \pmod{6}
\end{equation*}
Parenthesizing with $2^{2k+1}$ yields
\begin{equation*}
    2^{2k+1} + 2^{2k+1}(2^2-1) \equiv 2 \pmod{6}
\end{equation*}
which is
\begin{equation*}
    2^{2k+1} + 2^{2k}\times 6 \equiv 2 \pmod{6}
\end{equation*}
$6 \mid 2^{2k}\times 6$ and $2^{2k+1} \equiv 2 \pmod{6}$ is the congruence for $x=k$ therefore this lemma is proven by induction.
\end{proof}

\begin{lemma}
\label{LEMMA: 2^2x+2 equiv 4 mod 6}
$2^{2x+2} \equiv 4 \pmod{6}$ where $x \in \mathbb{N}_0$.
\end{lemma}
\begin{proof}
At $x=0$ we get $4 \equiv 4 \pmod{6}$. At $x=1$ we get $16 \equiv 4 \pmod{6}$. Assume that for $x=k$ the congruence
\begin{equation*}
    2^{2k+2} \equiv 4 \pmod{6}
\end{equation*}
is correct. Looking at $x=k+1$ we have
\begin{equation*}
    2^{2(k+1)+2} \equiv 4 \pmod{6}
\end{equation*}
Adding $2^{2k+2}-2^{2k+2}$ to this expression yields
\begin{equation*}
    2^{2k+2+2} + 2^{2k+2} - 2^{2k+2} \equiv 4 \pmod{6}
\end{equation*}
Parenthesizing with $2^{2k+2}$ yields
\begin{equation*}
    2^{2k+2} - 2^{2k+2}(2^2-1) \equiv 4 \pmod{6}
\end{equation*}
which is
\begin{equation*}
    2^{2k+2} + 2^{2k+1}\times 6 \equiv 4 \pmod{6}
\end{equation*}
$6 \mid 2^{2k+1}\times 6$ and $2^{2k+2} \equiv 4 \pmod{6}$ is the congruence for $x=k$ therefore this lemma is proven by induction.
\end{proof}

\begin{lemma}
\label{LEMMA: (2^6x+2 - 1) / 3 equiv 1 mod 6}
\begin{equation}
    \frac{2^{6x+2}-1}{3} \equiv 1 \pmod{6}
\end{equation}
 where $x \in \mathbb{N}_0$.
\end{lemma}
\begin{proof}
At $x=0$ we get $1 \equiv 1 \pmod{6}$. At $x=1$ we get $85 \equiv 1 \pmod{6}$. Assume that for $x=k$ the congruence
\begin{equation*}
    \frac{2^{6k+2}-1}{3} \equiv 1 \pmod{6}
\end{equation*}
is correct. Looking at $x=k+1$ we have
\begin{equation*}
    \frac{2^{6(k+1)+2}-1}{3} \equiv 1 \pmod{6}
\end{equation*}
Adding $2^{6k+2}-2^{6k+2}$ to $2^{6(k+1)+2}$ yields
\begin{equation*}
    \frac{2^{6(k+1)+2}+2^{6k+2}-2^{6k+2}-1}{3} \equiv 1 \pmod{6}
\end{equation*}
Parenthesizing with $2^{6k+2}$ yields
\begin{equation*}
    \frac{2^{6k+2}(2^6-1)}{3} + \frac{2^{6k+2}-1}{3} \equiv 1 \pmod{6}
\end{equation*}
which is
\begin{equation*}
    2^{6k+2}\times 42 + \frac{2^{6k+2}-1}{3} \equiv 1 \pmod{6}
\end{equation*}
$6 \mid 2^{6k+2}\times 42$ and $\frac{2^{6k+2}-1}{3} \equiv 1 \pmod{6}$ is the congruence for $x=k$ therefore this lemma is proven by induction.
\end{proof}

\begin{lemma}
\label{LEMMA: (2^6x+4 - 1) / 3 equiv 5 mod 6}
\begin{equation}
    \frac{2^{6x+4}-1}{3} \equiv 5 \pmod{6}
\end{equation}
\end{lemma}
\begin{proof}
Notice that this is equal to
\begin{equation*}
    \frac{4\times 2^{6x+2} - 1}{3} \equiv 5 \pmod{6}
\end{equation*}
therefore  
\begin{equation*}
    \frac{2^{6x+2} - 1}{3} + 2^{6x+2}  \equiv 5 \pmod{6}
\end{equation*}
According to lemma \ref{LEMMA: (2^6x+2 - 1) / 3 equiv 1 mod 6} we can say $\frac{2^{2(3x)+2} - 1}{3} \equiv 1 \pmod{6}$ and according to lemma \ref{LEMMA: 2^2x+2 equiv 4 mod 6} we can say
$2^{2(3x)+2} \equiv 4 \pmod{6}$ therefore $\frac{2^{6x+2} - 1}{3} + 2^{6x+2}  \equiv 1+4 \pmod{6}$.
\end{proof}

\begin{lemma}
\label{LEMMA: (2^6x+6 - 1) / 3 equiv 3 mod 6}
\begin{equation}
    \frac{2^{6x+6}-1}{3} \equiv 3 \pmod{6}
\end{equation}
\end{lemma}
\begin{proof}
Notice that this is equal to
\begin{equation*}
    \frac{4\times 2^{6x+4} - 1}{3} \equiv 3 \pmod{6}
\end{equation*}
therefore  
\begin{equation*}
    \frac{2^{6x+4} - 1}{3} + 2^{6x+4}  \equiv 3 \pmod{6}
\end{equation*}
According to lemma \ref{LEMMA: (2^6x+4 - 1) / 3 equiv 5 mod 6} we can say $\frac{2^{6x+4} - 1}{3} \equiv 5 \pmod{6}$ and according to lemma \ref{LEMMA: 2^2x+2 equiv 4 mod 6} we can say
$2^{2(3x+1)+2} \equiv 4 \pmod{6}$ therefore $\frac{2^{6x+4} - 1}{3} + 2^{6x+4}  \equiv 5+4 \pmod{6}$.
\end{proof}

\begin{theorem}
\label{THEO: RR_1}
Let $a, b \in \mathbb{N}_0$ and $c \in \{1, 3, 5\}$. The function $RR_1 : \mathbb{N}_0 \times \mathbb{N}_0 \times \{1, 3, 5\} \xrightarrow{} \mathbb{O}_{R1}$ is defined as
\begin{equation}
\label{EQ: RR_1}
    RR_1(a,b,c)=\frac{(6a+1)\times2^{2(3b+(c-1+(a \bmod{3}))\bmod{3})+2}-1}{3}
\end{equation}
where $RR_1(a,b,c) \equiv c \pmod{6}$.
\end{theorem}
\begin{proof}
We will consider the cases of $a$ in modulo 3 and $c$. $a \bmod 3$ can be $0, 1$ and $2$. We already defined $c \in \{1, 3, 5\}$. In total we have 9 pairs of $(a,c)$ which are $(0,1), (0,3), (0,5), (1,1), (1,3), (1,5), (2,1), (2,3), (2,5)$. We will look at equation \eqref{EQ: RR_1} using these pairs and look at the congruence $RR_1(a,b,c) \equiv c \pmod{6}$. Let us set $c=1$ first to show $RR_1(a,b,1) \equiv 1 \pmod{6}$.
\begin{equation*}
    RR_1(a,b,1) = \frac{(6a+1)2^{2(3b+((a \bmod{3}))\bmod{3})+2}-1}{3}
\end{equation*}
Let $a = 3k, k\in\mathbb{N}_0$ so $a \equiv 0 \pmod{3}$,
\begin{flalign*}
    RR_1(3k,b,1) &= \frac{(6(3k)+1)2^{2(3b+((3k \bmod{3}))\bmod{3})+2}-1}{3} \\
     RR_1(3k,b,1) &= \frac{(18k+1)2^{2(3b)+2}-1}{3} \\
      RR_1(3k,b,1) &= 6k\times2^{2(3b)+2} + \frac{2^{6b+2}-1}{3}
\end{flalign*}
$6 \mid 6k\times2^{2(3b)+2}$ and lemma \ref{LEMMA: (2^6x+2 - 1) / 3 equiv 1 mod 6} shows that $\frac{2^{6b+2}-1}{3} \equiv 1 \pmod{3}$, therefore $RR_1(3k, b, 1) \equiv 1 \pmod{6}$.
Let $a = 3k+1, k\in\mathbb{N}_0$ so $a \equiv 1 \pmod{3}$,
\begin{flalign*}
    RR_1(3k+1,b,1) &= \frac{(6(3k+1)+1)2^{2(3b+(((3k+1) \bmod{3}))\bmod{3})+2}-1}{3} \\
     RR_1(3k+1,b,1) &= \frac{(18k+6+1)2^{2(3b+1)+2}-1}{3} \\
      RR_1(3k+1,b,1) &= 6k\times2^{2(3b+1)+2} + 2\times2^{2(3b+1)+2} + \frac{2^{2(3b+1)+2}-1}{3} \\
      RR_1(3k+1,b,1) &= 6k\times2^{2(3b+1)+2} + 2^{2(3b+2)+1} + \frac{2^{6b+4}-1}{3}
\end{flalign*}
$6 \mid 6k\times2^{2(3b)+2}$, lemma \ref{LEMMA: 2^2x+1 equiv 2 mod 6} shows that $ 2^{2(3b+2)+1} \equiv 2 \pmod{6}$ and lemma \ref{LEMMA: (2^6x+4 - 1) / 3 equiv 5 mod 6} shows that $\frac{2^{6b+4}-1}{3} \equiv 5 \pmod{6}$, therefore $RR_1(3k+1, b, 1) \equiv 2+5 \pmod{6}$ so $RR_1(3k+1, b, 1) \equiv 1 \pmod{6}$.
Let $a = 3k+2, k\in\mathbb{N}_0$ so $a \equiv 2 \pmod{3}$,
\begin{flalign*}
    RR_1(3k+2,b,1) &= \frac{(6(3k+2)+1)2^{2(3b+(((3k+2) \bmod{3}))\bmod{3})+2}-1}{3} \\
     RR_1(3k+2,b,1) &= \frac{(18k+12+1)2^{2(3b+2)+2}-1}{3} \\
      RR_1(3k+2,b,1) &= 6k\times2^{2(3b+2)+2} + 4\times2^{2(3b+2)+2} + \frac{2^{2(3b+2)+2}-1}{3} \\
      RR_1(3k+2,b,1) &= 6k\times2^{2(3b+2)+2} + 2^{2(3b+3)+2} + \frac{2^{6b+6}-1}{3}
\end{flalign*}
$6 \mid 6k\times2^{2(3b+2)+2}$, lemma \ref{LEMMA: 2^2x+2 equiv 4 mod 6} shows that $ 2^{2(3b+3)+2} \equiv 4 \pmod{6}$ and lemma \ref{LEMMA: (2^6x+6 - 1) / 3 equiv 3 mod 6} shows that $\frac{2^{6b+6}-1}{3} \equiv 3 \pmod{6}$, therefore $RR_1(3k+2, b, 1) \equiv 4+3 \pmod{6}$ so $RR_1(3k+2, b, 1) \equiv 1 \pmod{6}$. This proves that $RR_1(a,b,1) \equiv 1 \pmod{6}$. Next, we look at $c=3$ to show $RR_1(a,b,3) \equiv 3 \pmod{6}$.
\begin{equation*}
    RR_1(a,b,3) = \frac{(6a+1)2^{2(3b+(2+(a \bmod{3}))\bmod{3})+2}-1}{3}
\end{equation*}
Let $a = 3k, k\in\mathbb{N}_0$ so $a \equiv 0 \pmod{3}$,
\begin{flalign*}
    RR_1(3k,b,3) &= \frac{(6(3k)+1)2^{2(3b+(2+(3k \bmod{3}))\bmod{3})+2}-1}{3} \\
     RR_1(3k,b,3) &= \frac{(18k+1)2^{2(3b+2)+2}-1}{3} \\
      RR_1(3k,b,3) &= 6k\times2^{2(3b+2)+2} + \frac{2^{6b+6}-1}{3}
\end{flalign*}
$6 \mid  6k\times2^{2(3b+2)+2}$ and lemma \ref{LEMMA: (2^6x+6 - 1) / 3 equiv 3 mod 6} shows that $\frac{2^{6b+6}-1}{3} \equiv 3 \pmod{6}$ therefore $RR_1(3k,b,3) \equiv 3 \pmod{6}$. Let $a = 3k+1, k\in\mathbb{N}_0$ so $a \equiv 1 \pmod{3}$,
\begin{flalign*}
    RR_1(3k+1,b,3) &= \frac{(6(3k+1)+1)2^{2(3b+(2+((3k+1) \bmod{3}))\bmod{3})+2}-1}{3} \\
     RR_1(3k+1,b,3) &= \frac{(18k+6+1)2^{2(3b)+2}-1}{3} \\
      RR_1(3k+1,b,3) &= 6k\times2^{2(3b)+2} + 2\times2^{2(3b)+2} + \frac{2^{2(3b)+2}-1}{3} \\
      RR_1(3k+1,b,3) &= 6k\times2^{2(3b)+2} + 2^{2(3b+1)+1} + \frac{2^{6b+2}-1}{3}
\end{flalign*}
$6 \mid 6k\times2^{2(3b)+2}$, lemma \ref{LEMMA: 2^2x+1 equiv 2 mod 6} shows that $2^{2(3b+1)+1} \equiv 2 \pmod{6}$ and lemma \ref{LEMMA: (2^6x+2 - 1) / 3 equiv 1 mod 6} shows that $\frac{2^{6b+2}-1}{3} \equiv 1 \pmod{6}$ therefore $RR_1(3k+1,b,3) \equiv 1+2 \pmod{6}$ so $RR_1(3k+1,b,3) \equiv 3 \pmod{6}$. Let $a = 3k+2, k\in\mathbb{N}_0$ so $a \equiv 2 \pmod{3}$,
\begin{flalign*}
    RR_1(3k+2,b,3) &= \frac{(6(3k+2)+1)2^{2(3b+(2+((3k+2) \bmod{3}))\bmod{3})+2}-1}{3} \\
     RR_1(3k+2,b,3) &= \frac{(18k+12+1)2^{2(3b+1)+2}-1}{3} \\
      RR_1(3k+2,b,3) &= 6k\times2^{2(3b+1)+2} + 4\times2^{2(3b+1)+2} + \frac{2^{2(3b+1)+2}-1}{3} \\
      RR_1(3k+2,b,3) &= 6k\times2^{2(3b+1)+2} + 2^{2(3b+2)+2} + \frac{2^{6b+4}-1}{3}
\end{flalign*}
$6 \mid 6k\times2^{2(3b+1)+2}$, lemma \ref{LEMMA: 2^2x+2 equiv 4 mod 6} shows that $2^{2(3b+2)+2} \equiv 4 \pmod{6}$ and lemma \ref{LEMMA: (2^6x+4 - 1) / 3 equiv 5 mod 6} shows that $\frac{2^{6b+4}-1}{3} \equiv 5 \pmod{6}$ therefore $RR_1(3k+2,b,3) \equiv 4 + 5 \pmod{6}$ so $RR_1(3k+2,b,3) \equiv 3 \pmod{6}$. This proves that $RR_1(a,b,3) \equiv 3 \pmod{6}$. Finally, we look at $c=5$ to show $RR_1(a,b,5) \equiv 5 \pmod{6}$.
\begin{equation*}
    RR_1(a,b,5)=\frac{(6a+1)\times2^{2(3b+(4+(a \bmod{3}))\bmod{3})+2}-1}{3}
\end{equation*}
Let $a = 3k, k\in\mathbb{N}_0$ so $a \equiv 0 \pmod{3}$,
\begin{flalign*}
    RR_1(3k,b,5) &= \frac{(6(3k)+1)2^{2(3b+(4+(3k \bmod{3}))\bmod{3})+2}-1}{3} \\
     RR_1(3k,b,5) &= \frac{(18k+1)2^{2(3b+1)+2}-1}{3} \\
      RR_1(3k,b,5) &= 6k\times2^{2(3b+1)+2} + \frac{2^{6b+4}-1}{3}
\end{flalign*}
$6 \mid 6k\times2^{2(3b+1)+2}$ and lemma \ref{LEMMA: (2^6x+4 - 1) / 3 equiv 5 mod 6} shows that $\frac{2^{6b+4}-1}{3} \equiv 5 \pmod{6}$ therefore $RR_1(3k,b,5) \equiv 5 \pmod{6}$. Let $a = 3k+1, k\in\mathbb{N}_0$ so $a \equiv 1 \pmod{3}$,
\begin{flalign*}
    RR_1(3k+1,b,5) &= \frac{(6(3k+1)+1)2^{2(3b+(4+((3k+1) \bmod{3}))\bmod{3})+2}-1}{3} \\
     RR_1(3k+1,b,5) &= \frac{(18k+6+1)2^{2(3b+2)+2}-1}{3} \\
      RR_1(3k+1,b,5) &= 6k\times2^{2(3b+2)+2} + 2\times2^{2(3b+2)+2} + \frac{2^{2(3b+2)+2}-1}{3} \\
      RR_1(3k+1,b,5) &= 6k\times2^{2(3b+2)+2} + 2^{2(3b+3)+1} + \frac{2^{6b+6}-1}{3}
\end{flalign*}
$6 \mid 6k\times2^{2(3b+2)+2}$, lemma \ref{LEMMA: 2^2x+1 equiv 2 mod 6} shows that $2^{2(3b+3)+1} \equiv 2 \pmod{6}$ and lemma \ref{LEMMA: (2^6x+6 - 1) / 3 equiv 3 mod 6} shows that $\frac{2^{6b+6}-1}{3} \equiv 3 \pmod{6}$ therefore $RR_1(3k+1,b,5) \equiv 2 + 3 \pmod{6}$ so $RR_1(3k+1,b,5) \equiv 5 \pmod{6}$. Let $a = 3k+2, k\in\mathbb{N}_0$ so $a \equiv 2 \pmod{3}$,
\begin{flalign*}
    RR_1(3k+2,b,5) &= \frac{(6(3k+2)+1)2^{2(3b+(4+((3k+2) \bmod{3}))\bmod{3})+2}-1}{3} \\
     RR_1(3k+2,b,5) &= \frac{(18k+12+1)2^{2(3b)+2}-1}{3} \\
      RR_1(3k+2,b,5) &= 6k\times2^{2(3b)+2} + 4\times2^{2(3b)+2} + \frac{2^{2(3b)+2}-1}{3} \\
      RR_1(3k+2,b,5) &= 6k\times2^{2(3b)+2} + 2^{2(3b+1)+2} + \frac{2^{6b+2}-1}{3}
\end{flalign*}
$6 \mid 6k\times2^{2(3b)+2}$, lemma \ref{LEMMA: 2^2x+2 equiv 4 mod 6} shows that $2^{2(3b+1)+2} \equiv 4 \pmod{6}$ and lemma \ref{LEMMA: (2^6x+2 - 1) / 3 equiv 1 mod 6} shows that $\frac{2^{6b+2}-1}{3} \equiv 1 \pmod{6}$ therefore $RR_1(3k+2,b,5) \equiv 4+1$ so $RR_1(3k+2,b,5) \equiv 5 \pmod{6}$. This proves that $RR_1(a,b,5) \equiv 5 \pmod{6}$. Having proved $RR_1(a,b,1) \equiv 1 \pmod{6}$ and $RR_1(a,b,3) \equiv 3 \pmod{6}$ too, we are able to say $RR_1(a,b,c) \equiv c \pmod{6}$ where $c \in \{1, 3, 5\}$. We should also show that the range of function $RR_1$ is $\mathbb{O}_{R1}$. Comparing $RR_1$ to $R_1$, only difference is the exponent of 2. The exponent in $RR_1$ is $2(3b+(c-1+(a \bmod 3)) \bmod 3)+2$ and in $R_1$ it is $2b+2$. For $R_1$ we know $b \in \mathbb{N}_0$, therefore we have to show $3b+(c-1+(a \bmod 3)) \bmod 3$ can produce all numbers in $\mathbb{N}_0$. This can be seen if we look at the proofs for all 9 pairs of $(a,c)$ above. As we are working with the pairs, we see that the exponent becomes $3b, 3b+1$ and $3b+2$. Since $b \in \mathbb{N}_0$ it is clear $3b, 3b+1$ and $3b+2$ will produce numbers in $\mathbb{N}_0$.
\end{proof}

\begin{theorem}
\label{THEO: RR_5}
Let $a, b \in \mathbb{N}_0$ and $c \in \{1, 3, 5\}$. The function $RR_5 : \mathbb{N}_0 \times \mathbb{N}_0 \times \{1, 3, 5\} \xrightarrow{} \mathbb{O}_{R5}$ is defined as
\begin{equation}
\label{EQ: RR_5}
    RR_5(a,b,c)=\frac{(6a+5)\times2^{2(3b+(c-(a \bmod{3}))\bmod{3})+1}-1}{3}
\end{equation}
where $RR_5(a,b,c) \equiv c \pmod{6}$.
\end{theorem}
\begin{proof}
The proof will be similar to that of theorem \ref{THEO: RR_1}. Let us set $c=1$ first to show $RR_5(a,b,1) \equiv 1 \pmod{6}$. 
\begin{equation*}
    RR_5(a,b,1) = \frac{(6a+5)2^{2(3b+(1-(a \bmod{3}))\bmod{3})+1}-1}{3}
\end{equation*}
Let $a = 3k, k \in \mathbb{N}_0$ so $a \equiv 0 \pmod{3}$.
\begin{flalign*}
    RR_5(3k,b,1) &= \frac{(6(3k)+5)2^{2(3b+(1-(3k \bmod{3}))\bmod{3})+1}-1}{3} \\
    RR_5(3k,b,1) &= \frac{(6(3k)+5)2^{2(3b+1)+1}-1}{3} \\
    RR_5(3k,b,1) &= \frac{(18k+3+2)2^{2(3b+1)+1}-1}{3} \\
    RR_5(3k,b,1) &= 6k \times 2^{2(3b+1)+1} + 2^{2(3b+1)+1} +\frac{2 \times 2^{6b+3}-1}{3} \\
    RR_5(3k,b,1) &= 6k \times 2^{2(3b+1)+1} + 2^{2(3b+1)+1} +\frac{ 2^{6b+4}-1}{3}
\end{flalign*}
$6 \mid6k \times 2^{2(3b+1)+1}$, lemma \ref{LEMMA: 2^2x+1 equiv 2 mod 6} shows that $2^{2(3b+1)+1} \equiv 2 \pmod{6}$ and lemma \ref{LEMMA: (2^6x+4 - 1) / 3 equiv 5 mod 6} shows that $\frac{ 2^{6b+4}-1}{3} \equiv 5 \pmod{6}$ therefore $RR_5(3k,b,1) \equiv 2+5 \pmod{6}$ so $RR_5(3k,b,1) \equiv 1 \pmod{6}$. Let $a = 3k+1, k \in \mathbb{N}_0$ so $a \equiv 1 \pmod{3}$.
\begin{flalign*}
    RR_5(3k+1,b,1) &= \frac{(6(3k+1)+5)2^{2(3b+(1-((3k+1) \bmod{3}))\bmod{3})+1}-1}{3} \\
    RR_5(3k+1,b,1) &= \frac{(6(3k+1)+5)2^{2(3b)+1}-1}{3} \\
    RR_5(3k+1,b,1) &= \frac{(18k+9+2)2^{2(3b)+1}-1}{3} \\
    RR_5(3k+1,b,1) &= 6k \times 2^{2(3b)+1} + 3\times2^{2(3b)+1}+\frac{2 \times 2^{6b+1}-1}{3} \\
    RR_5(3k+1,b,1) &= 6k \times 2^{2(3b)+1} + 6\times2^{2(3b)} +\frac{ 2^{6b+2}-1}{3}
\end{flalign*}
$6 \mid 6k \times 2^{2(3b)+1} + 6\times2^{2(3b)}$ and lemma \ref{LEMMA: (2^6x+2 - 1) / 3 equiv 1 mod 6} shows that $\frac{ 2^{6b+2}-1}{3} \equiv 1 \pmod{6}$ therefore $RR_5(3k+1,b,1) \equiv 1 \pmod{6}$. Let $a = 3k+2, k \in \mathbb{N}_0$ so $a \equiv 2 \pmod{3}$. 
\begin{flalign*}
    RR_5(3k+2,b,1) &= \frac{(6(3k+2)+5)2^{2(3b+(1-((3k+2) \bmod{3}))\bmod{3})+1}-1}{3} \\
    RR_5(3k+2,b,1) &= \frac{(6(3k+2)+5)2^{2(3b+2)+1}-1}{3} \\
    RR_5(3k+2,b,1) &= \frac{(18k+12+3+2)2^{2(3b+2)+1}-1}{3} \\
    RR_5(3k+2,b,1) &= 6k \times 2^{2(3b+2)+1} + 4\times2^{2(3b+2)+1}+2^{2(3b+2)+1}+\frac{2 \times 2^{6b+5}-1}{3} \\
    RR_5(3k+2,b,1) &= 6k \times 2^{2(3b+2)+1} + 2^{2(3b+3)+1}+2^{2(3b+2)+1}+\frac{ 2^{6b+6}-1}{3}
\end{flalign*}
$6 \mid 6k \times 2^{6b+5}$, lemma \ref{LEMMA: 2^2x+1 equiv 2 mod 6} shows that $2^{2(3b+3)+1}+2^{2(3b+2)+1} \equiv 2 + 2 \pmod{6}$ and lemma \ref{LEMMA: (2^6x+6 - 1) / 3 equiv 3 mod 6} shows $\frac{ 2^{6b+6}-1}{3} \equiv 3 \pmod{6}$ therefore $RR_5(3k+2,b,1) \equiv 2 + 2 + 3 \pmod{6}$ so $RR_5(3k+2,b,1) \equiv 1 \pmod{6}$. This proves that $RR_5(a,b,1) \equiv 1 \pmod{6}$. Next, we look at $c=3$ to show $RR_5(a,b,3) \equiv 3 \pmod{6}$.
\begin{equation*}
    RR_5(a,b,3) = \frac{(6a+5)2^{2(3b+(3-(a \bmod{3}))\bmod{3})+1}-1}{3}
\end{equation*}
Let $a = 3k, k \in \mathbb{N}_0$ so $a \equiv 0 \pmod{3}$. 
\begin{flalign*}
    RR_5(3k,b,3) &= \frac{(6(3k)+5)2^{2(3b+(3-(3k \bmod{3}))\bmod{3})+1}-1}{3} \\
    RR_5(3k,b,3) &= \frac{(18k+3+2)2^{2(3b)+1}-1}{3} \\
    RR_5(3k,b,3) &= 6k \times 2^{2(3b)+1} + 2^{2(3b)+1} + \frac{2\times 2^{6b+1} -1}{3} \\
    RR_5(3k,b,3) &= 6k \times 2^{2(3b)+1} + 2^{2(3b)+1} + \frac{2^{6b+2} -1}{3}
\end{flalign*}
$6 \mid 6k \times 2^{2(3b)+1}$, lemma \ref{LEMMA: 2^2x+1 equiv 2 mod 6} shows that $ 2^{2(3b)+1} \equiv 2 \pmod{6}$ and lemma \ref{LEMMA: (2^6x+2 - 1) / 3 equiv 1 mod 6} shows $\frac{2^{6b+2} -1}{3} \equiv 1 \pmod{6}$ therefore $RR_5(3k,b,3) \equiv 2 + 1 \pmod{6}$ so $RR_5(3k,b,3) \equiv 3 \pmod{6}$. Let $a = 3k+1, k \in \mathbb{N}_0$ so $a \equiv 1 \pmod{3}$. 
\begin{flalign*}
    RR_5(3k+1,b,3) &= \frac{(6(3k+1)+5)2^{2(3b+(3-((3k+1) \bmod{3}))\bmod{3})+1}-1}{3} \\
    RR_5(3k+1,b,3) &= \frac{(18k+9+2)2^{2(3b+2)+1}-1}{3} \\
    RR_5(3k+1,b,3) &= 6k \times 2^{2(3b+2)+1} + 3\times 2^{2(3b+2)+1} + \frac{2\times 2^{6b+5} -1}{3} \\
    RR_5(3k+1,b,3) &= 6k \times 2^{2(3b+2)+1} + 6\times 2^{2(3b+2)} + \frac{2^{6b+6} -1}{3}
\end{flalign*}
$6 \mid 6k \times 2^{2(3b+2)+1} + 6\times 2^{2(3b+2)}$ and lemma \ref{LEMMA: (2^6x+6 - 1) / 3 equiv 3 mod 6} shows that $\frac{2^{6b+6} -1}{3} \equiv 3 \pmod{6}$ therefore $RR_5(3k+1,b,3) \equiv 3 \pmod{6}$. Let $a = 3k+2, k \in \mathbb{N}_0$ so $a \equiv 2 \pmod{3}$.
\begin{flalign*}
    RR_5(3k+2,b,3) &= \frac{(6(3k+2)+5)2^{2(3b+(3-((3k+2) \bmod{3}))\bmod{3})+1}-1}{3} \\
    RR_5(3k+2,b,3) &= \frac{(18k+12+3+2)2^{2(3b+1)+1}-1}{3} \\
    RR_5(3k+2,b,3) &= 6k \times 2^{2(3b+1)+1} + 4\times 2^{2(3b+1)+1} + 2^{2(3b+1)+1} + \frac{2\times 2^{6b+3} -1}{3} \\
    RR_5(3k+2,b,3) &= 6k \times 2^{2(3b+1)+1} + 2^{2(3b+2)+1} + 2^{2(3b+1)+1} + \frac{ 2^{6b+4} -1}{3}
\end{flalign*}
$6 \mid 6k \times 2^{2(3b+1)+1}$, lemma \ref{LEMMA: 2^2x+1 equiv 2 mod 6} shows that $2^{2(3b+2)+1} + 2^{2(3b+1)+1} \equiv 2 + 2 \pmod{6}$ and lemma \ref{LEMMA: (2^6x+4 - 1) / 3 equiv 5 mod 6} shows that $\frac{ 2^{6b+4} -1}{3} \equiv 5 \pmod{6}$ therefore $RR_5(3k+2,b,3) \equiv 2 + 2 + 5 \pmod{6}$ so $RR_5(3k+2,b,3) \equiv 3 \pmod{6}$. This proves that $RR_5(a,b,3) \equiv 3 \pmod{6}$. Finally we set $c=5$ to show $RR_5(a,b,5) \equiv 5 \pmod{6}$. 
\begin{equation*}
    RR_5(a,b,5) = \frac{(6a+5)2^{2(3b+(5-(a \bmod{3}))\bmod{3})+1}-1}{3}
\end{equation*}
Let $a = 3k, k \in \mathbb{N}_0$ so $a \equiv 0 \pmod{3}$. 
\begin{flalign*}
    RR_5(3k,b,5) &= \frac{(6(3k)+5)2^{2(3b+(5-(3k \bmod{3}))\bmod{3})+1}-1}{3} \\
    RR_5(3k,b,5) &= \frac{(18k+3+2)2^{2(3b+2)+1}-1}{3} \\
    RR_5(3k,b,5) &= 6k\times 2^{2(3b+2)+1} + 2^{2(3b+2)+1} + \frac{2\times2^{6b+5} - 1}{3} \\
    RR_5(3k,b,5) &= 6k\times 2^{2(3b+2)+1} + 2^{2(3b+2)+1} + \frac{2^{6b+6} - 1}{3} 
\end{flalign*}
$6 \mid  6k\times 2^{2(3b+2)+1} $, lemma \ref{LEMMA: 2^2x+1 equiv 2 mod 6} shows that $2^{2(3b+2)+1} \equiv 2 \pmod{6}$ and lemma \ref{LEMMA: (2^6x+6 - 1) / 3 equiv 3 mod 6} shows $\frac{2^{6b+6} - 1}{3} \equiv 3 \pmod{6}$ therefore $RR_5(3k,b,5) \equiv 2 + 3 \pmod{6}$ so $RR_5(3k,b,5) \equiv 5 \pmod{6}$. Let $a = 3k+1, k \in \mathbb{N}_0$ so $a \equiv 1 \pmod{3}$. 
\begin{flalign*}
    RR_5(3k+1,b,5) &= \frac{(6(3k+1)+5)2^{2(3b+(5-((3k+1) \bmod{3}))\bmod{3})+1}-1}{3} \\
    RR_5(3k+1,b,5) &= \frac{(6(3k+1)+5)2^{2(3b+1)+1}-1}{3} \\
    RR_5(3k+1,b,5) &= \frac{(18k+9+2)2^{2(3b+1)+1}-1}{3} \\
    RR_5(3k+1,b,5) &= 6k\times2^{2(3b+1)+1} + 3\times 2^{2(3b+1)+1} + \frac{2\times 2^{6b+3}}{3}  \\
     RR_5(3k+1,b,5) &= 6k\times2^{2(3b+1)+1} + 6\times 2^{2(3b+1)} + \frac{2^{6b+4}}{3}  \\
\end{flalign*}
$6 \mid 6k\times2^{2(3b+1)+1} + 6\times 2^{2(3b+1)} $ and lemma \ref{LEMMA: (2^6x+4 - 1) / 3 equiv 5 mod 6} shows $\frac{2^{6b+4}}{3} \equiv 5 \pmod{6}$ therefore $RR_5(3k+1,b,5) \equiv 5 \pmod{6}$. Let $a = 3k+2, k \in \mathbb{N}_0$ so $a \equiv 2 \pmod{3}$. 
\begin{flalign*}
    RR_5(3k+2,b,5) &= \frac{(6(3k+2)+5)2^{2(3b+(5-((3k+2) \bmod{3}))\bmod{3})+1}-1}{3} \\
    RR_5(3k+2,b,5) &= \frac{(18k+12+3+2)2^{2(3b)+1}-1}{3} \\
    RR_5(3k+2,b,5) &= 6k\times 2^{2(3b)+1} + 4\times 2^{2(3b)+1} + 2^{2(3b)+1} + \frac{2\times 2^{6b+1}-1}{3} \\
    RR_5(3k+2,b,5) &= 6k\times 2^{2(3b)+1} + 2^{2(3b+1)+1} + 2^{2(3b)+1} + \frac{2^{6b+2}-1}{3}
\end{flalign*}
$6 \mid 6k\times 2^{2(3b)+1}$, lemma \ref{LEMMA: 2^2x+1 equiv 2 mod 6} shows that $2^{2(3b+1)+1} + 2^{2(3b)+1} \equiv 2 + 2 \pmod{6}$ and lemma \ref{LEMMA: (2^6x+2 - 1) / 3 equiv 1 mod 6} shows that $\frac{2^{6b+2}-1}{3} \equiv 1 \pmod{6}$ therefore $RR_5(3k+2,b,5) \equiv 2+2+1 \pmod{6}$ so $RR_5(3k+2,b,5) \equiv 5 \pmod{6}$. This proves that $RR_5(a,b,5) \equiv 5 \pmod{6}$. Having proved $RR_5(a,b,1) \equiv 1 \pmod{6}$ and $RR_5(a,b,3) \equiv 3 \pmod{6}$ too, we are able to say $RR_5(a,b,c) \equiv c \pmod{6}$ where $c \in \{1, 3, 5\}$. We should also show that the range of function $RR_5$ is $\mathbb{O}_{R5}$. Comparing $RR_5$ to $R_5$the only difference is the exponent of $2$. The exponent in $RR_5$ is $2(3b+(c-(a \bmod{3}))\bmod{3})+1$ and in $R_5$ it is $2b+1$. For $R_5$ we know $b \in \mathbb{N}_0$, therefore we have to show $2(3b+(c-(a \bmod{3}))\bmod{3})+1$ can produce all numbers in $\mathbb{N}_0$. This can be seen if we look at the proofs for all 9 pairs of $(a,c)$ above. As we are working with the pairs, we see that the exponent becomes $3b$, $3b+1$ and $3b+2$. Since $b \in \mathbb{N}_0$ it is clear $3b, 3b+1$ and $3b+2$ will produce numbers in $\mathbb{N}_0$.
\end{proof}

%% file: Sections/4.2-OneTree.tex
\subsection{One-tree using reverse reduced Collatz function}
Consider a graph $G = (V, E)$ where $V = \mathbb{O}_+$  and $E = \{(n,m) : C(n) = m \text{ where } n,m \in V\}$. This graph is most notably known as Collatz Graph \cite{lagarias-general}. If we omit the self-loop $(1, 1)$ from the Collatz Graph, the conjecture can be reformulated in this context: ``Is Collatz Graph a tree with root vertex $1$''? The aforementioned tree is also named One-tree in the context of Collatz conjecture, or $(3x+1)$-tree \cite{lagarias-bib}. 
\subsubsection{Constructing the Collatz graph}
Let $s = R(d)$. By definition of $R$ we know that $s$ does not have a single value, since the exponent $x$ in $R$ can have many values (see \eqref{EQ: Reverse Formula}). Moreover, $s$ does not have a value when $d \equiv 3 \pmod{6}$ (see lemma \ref{LEMMA: 6a+3 ise x yok}). We will consider $s$ separately for $d \equiv 1 \pmod{6}$ and $d \equiv 5 \pmod{6}$. We can do this using $R_1$ and $R_5$ functions. We get $s = R_1(a,b), d = 6a+1$ or $s = R_5(a,b), d = 6a+5$.
Suppose we calculate $s$ and $d$ by calculating $R_1(a,b)$ or $R_5(a,b)$. The order of calculation or the starting point does not matter, what matters is that we cover all possibilities for both functions where $a, b \in \mathbb{N}_0$. Suppose that we made a calculation and thereby have $s$ and $d$. Suppose we also have a graph $G = (V, E)$. We have 4 cases regarding $s$ and $d$ being an element of $V$:
\begin{enumerate}
    \item $s \not\in V, d \not\in V$. Neither $s$ or $d$ is present in our graph. What we do is, add $s$ and $d$ to our graph and connect $s$ to $d$. Our new graph is $G' = (V \cup \{s, d\}, E \cup \{(s, d)\})$. 
    \item $s \not\in V, d \in V$. This means that $d$ already exists in the graph, so we add $s$ to the graph and connect it to $d$. Our new graph is $G' = (V \cup \{s\}, E \cup \{(s, d)\})$.
    \item $s \in V, d \not\in V$. This means that $s$ already exists in the graph, so we add $d$ to the graph and connect it to $s$. Our new graph is $G' = (V \cup \{d\}, E \cup \{(s, d)\})$.
    \item $s \in V, d \in V$. This means that both numbers are present in the graph but they are not connected. Our new graph is $G' = (V, E \cup \{(s, d)\})$
\end{enumerate}
Suppose we use $G'$ instead of $G$ after every calculation of $s$ and $d$. Starting from an empty graph $G = (\varnothing, \varnothing)$ and calculating $s$, $d$ for all possibilities of $R_1(a,b)$ and $R_5(a,b)$ using the logic described above, we expect to get the Collatz graph. From theorem \ref{THEO: R1 U R5 = O} we can infer that a case where $s, d \in V$ and $(s, d) \in E$ at some calculation is not possible, because the theorem indicates that $s$ will always be unique, as for the graph construction in all four cases we add the edge $(s, d)$ to $E$, if $s$ is unique then the edge $(s, d)$ will be added to $E$ only once, for every unique occurrence of $s$, that is.
\begin{figure}[!ht]
    \centering
    \includegraphics[width=0.8\textwidth]{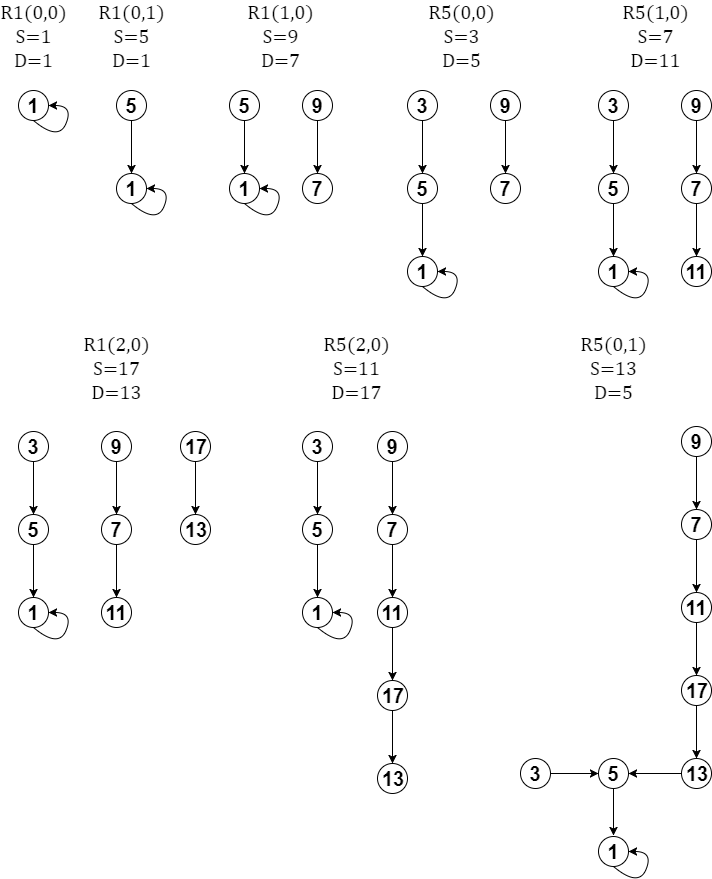}
    \caption{An example step-by-step construction of Collatz Tree, steps from left to right and top to bottom.}
    \label{FIG: Collatz Graph}
\end{figure}
\par A short example of the procedure is given in figure \ref{FIG: Collatz Graph}. There are 8 calculations, order being from left to right and top to bottom. First $R_1(0,0)$ is calculated which gives $s=1$, $d=1$. This is the only case where $s=d$, which is the trivial cycle $1\xrightarrow{2}1$. We add the vertex $1$ and make a self-loop. Next, we calculate $R_1(0,1)$ and get $s=5$, $d=1$. $1$ exists in the vertices (in fact it is the only one so far) but 5 does not. This is case 2 as described above. We add 5 and connect it to 1. The calculations go on like this; it is conjectured that in the end we will have a tree with root 1 (and the self-loop at 1 of course), similar to the graph shown at the bottom right corner of figure \ref{FIG: Collatz Graph}. 

%% file: Sections/5-Conclusion.tex
\section{Concluding Remarks}
We believe the fractional sum notations (see definitions \ref{DEF: FSN Definition} and \ref{DEF: IFSN Definition}) are a useful way of analyzing reduced Collatz trajectories. A continuation of the loops discussed in section \ref{SECT: Loops} could perhaps shed more light on the subject. The issue regarding equation \ref{EQ: Additive Formula} shows that it can perhaps be improved to resolve the problem, which would let us study Collatz sequences even better.

\par We have shown that the reverse reformulation can used to predetermine the numbers in modulo 6 (see theorems \ref{THEO: RR_1} and \ref{THEO: RR_5}). For example, one can avoid bumping into numbers $n$ such that $n \equiv 3 \pmod{6}$ (which are problematic for reverse reduced Collatz function, see \ref{LEMMA: 6a+3 ise x yok}) by using the modified functions \eqref{EQ: RR_1} and \eqref{EQ: RR_5}. 

\par For a trajectory $n \xrightarrow{a_1} C(n) \xrightarrow{a_2} \ldots \xrightarrow{a_j} 1$ we could use equation \eqref{EQ: Reverse Formula} to get $R(n) \xrightarrow{x} n \xrightarrow{a_1} C(n) \xrightarrow{a_2} \ldots \xrightarrow{a_j} 1$. We have shown how to use equation \eqref{EQ: Additive Formula} and equation \eqref{EQ: Additive Formula using Trivial Cycle} in theorem \ref{THEO: Additive Formula} to get $m \xrightarrow{a_1} C(m) \xrightarrow{a_2} \ldots \xrightarrow{a_j} C^j(m) \xrightarrow{a_{j+1}} 1$. This may provide an additional perspective on generating trajectories.


%% file: main.bbl
\begin{thebibliography}{10}

\bibitem{kerstin}
K.~{Andersson}.
\newblock {On the Boundedness of Collatz Sequences}.
\newblock {\em arXiv e-prints}, page arXiv:1403.7425, Mar 2014.

\bibitem{crandall}
R.~E. Crandall.
\newblock On the 3x+1 {P}roblem.
\newblock {\em Mathematics of Computation}, 32(144):1281--1292, October 1978.

\bibitem{ghelms}
G.~Helms.
\newblock Collatz-{I}ntro - {S}ome general remarks, 2004.
\newblock Accessed 12 March 2019.

\bibitem{lagarias-general}
J.~C. Lagarias.
\newblock The 3x+1 and {I}ts {G}eneralizations.
\newblock {\em American Mathematical Monthly}, 92(1):3--23, January 1985.

\bibitem{lagarias-bib}
J.~C. {Lagarias}.
\newblock {The 3x+1 Problem: An Annotated Bibliography, II (2000-2009)}.
\newblock {\em arXiv Mathematics e-prints}, page math/0608208, Aug 2006.

\bibitem{fcmotta}
F.~C. Motta, H.~R. de~Oliveira, and T.~A. Catalan.
\newblock An {A}nalysis of the {C}ollatz {C}onjecture.

\bibitem{ericroosendall}
E.~Roosendall.
\newblock On the 3x+1 {P}roblem, 2019.
\newblock Accessed 12 March 2019.

\bibitem{simons}
J.~L. Simons.
\newblock A simple (inductive) proof for the non-existence of 2-cycles of the
  3x+1 problem.
\newblock {\em Journal of Number Theory}, 123(1):10 -- 17, 2007.

\bibitem{steiner}
R.~P. Steiner.
\newblock A {T}heorem on the {S}yracuse {P}roblem.
\newblock In {\em A {T}heorem on the {S}yracuse {P}roblem}, pages 553--559. 7th
  Manitoba Conference on Numerical Mathematics, 1977.

\bibitem{ianstewart}
Ian Stewart.
\newblock {\em The {G}reat {M}athematical {P}roblems}, chapter~17, page 284.
\newblock Profile Books LTD, 2014.

\end{thebibliography}
